\documentstyle [12pt]{article}
 \let\tilde=\widetilde
\let\oldphi=\phi \def\e{{\rm e}} \def\phi{\varphi}
\def\eps{\varepsilon} \def\a{\alpha} \def\b{\beta}
 \def\G{\Gamma} \def\d{\delta}  
\def\l{\lambda}

 \def\res{\mathop{\rm res}\limits}
 \def\Sym{\mathop{\rm Sym}\limits}
 \def\dilog{\mathop{\rm Li}\nolimits_2}

 \def\dd{\partial}

\def\one#1{#1^{\raise5pt\hbox{$\scriptstyle\!\!\!\!1$}}\,{}}
\def\two#1{#1^{\raise5pt\hbox{$\scriptstyle\!\!\!\!2$}}\,{}}
\def\three#1{#1^{\raise5pt\hbox{$\scriptstyle\!\!\!\!3$}}\,{}}
 
 \def\abs#1{\left|#1\right|}
\def\id{\hbox{{1}\kern-.25em\hbox{\rm l}}} \def\intl{\int\limits}
 \def\comment#1{}
 
\def\endproof{\hfill\rule{2mm}{2mm}}
\def\beq{\begin{equation}} \def\eeq{\end{equation}}
\def\be{\begin{displaymath}} \def\ee{\end{displaymath}}
\def\bea{\begin{eqnarray}} \def\eea{\end{eqnarray}}
\def\beas{\begin{eqnarray*}} \def\eeas{\end{eqnarray*}}
\def\bds{\begin{description}} \def\eds{\end{description}}
\def\bmat{\left(\begin{array}} \def\emat{\end{array}\right)}
 \newtheorem{theo}{Theorem}
\newtheorem{prop}{Proposition} 
 \def\half{\frac{1}{2}}

\def\dfrac#1#2{\frac{\displaystyle #1}{\displaystyle #2}}
\def\Ref#1{(\ref{#1})}
\def\?{(?)\marginpar{|?}}

 \def\I{{\cal I}}  \def\M{{\cal M}}
  
\def\Lq#1#2#3{{\cal L}_q\!\left(#1;#2,#3\right)}

\renewcommand{\theequation}{\thesection.\arabic{equation}}
\newcounter{subequation}[equation]
\makeatletter
 
\expandafter\let\expandafter\reset@font\csname reset@font\endcsname
 
\def\subeqnarray{\arraycolsep1pt
    \def\@eqnnum\stepcounter##1{\stepcounter{subequation}%
        {\reset@font\rm(\theequation\alph{subequation})}}
\jot5mm     \eqnarray}

\makeatother\newcommand{\newsection}[1]{\vspace{10mm}
\pagebreak[3]\addtocounter{section}{1}\setcounter{equation}{0}
\setcounter{subsection}{0}\setcounter{footnote}{0}
 
\begin{flushleft}{\Large\bf \thesection. #1}
\end{flushleft}\nopagebreak\medskip\nopagebreak}
\newcommand{\newappendix}[1]{\vspace{10mm}\pagebreak[3]
\addtocounter{section}{1}\setcounter{equation}{0}
 
\begin{flushleft}{\Large\bf Appendix \thesection. #1}
\end{flushleft}\nopagebreak\medskip\nopagebreak}

\newfont{\bbd}{msbm10 scaled\magstep1} \def\C{\hbox{\bbd C}}                  
\def\R{\hbox{\bbd R}}  \def\Z{\hbox{\bbd Z}}

\newfont{\frak}{eufm10 scaled\magstep1}

\newcounter{remctr} \newenvironment{rem}
{\refstepcounter{remctr}%
{\bf Remark \arabic{remctr}\ }}{}
\def\sfrac#1#2{\frac{\scriptstyle#1}{\scriptstyle#2}}
\def\Symx{\Sym(x_1,x_2)}\def\Symy{\Sym(y_1,y_2)}
\def\Proof{{\bf Proof. }}

\newcommand{\tfrac}[2]{{\textstyle\frac{#1}{#2}}}
\begin{document}
\begin{flushright}
October 1996 \\
{\tt q-alg/9705006}
\end{flushright}
\begin{center} \LARGE\bf
Separation of variables and integral relations for special functions
\end{center}
\vskip 0.5cm
\begin{center}
V.B.~Kuznetsov
\footnote{On leave from: Department of Mathematical and Computational
Physics, Institute of Physics, St.~Petersburg University, St.~Petersburg
198904, Russian Federation.} \\
\vskip 0.2cm
Department of Applied Mathematical Studies,\\  
University of Leeds, Leeds LS2 9JT, UK\\
\vskip 0.7cm
E.K. Sklyanin \footnote{On leave from: Steklov Mathematical Institute,
Fontanka 27, St.~Petersburg 191011, Russian Federation.}
\vspace{.2cm} \\ 
Research Institute for Mathematical Sciences \\
Kyoto University, Kyoto 606, Japan
\end{center}
\vskip 1cm
\centerline{\it To the memory of Felix M. Arscott}
\vskip 1.5cm
\begin{center}
\bf Abstract
\end{center}
We show that the method of separation of variables gives a natural
generalisation of integral relations for classical special functions
of one variable. The approach is illustrated by giving a new
proof of the ``quadratic'' integral relations for the 
continuous $q$-ultraspherical
polynomials. The separating integral operator $M$ expressed in terms
of the Askey-Wilson operator is studied in detail: apart from writing
down the characteristic (``separation'') equations it satisfies, we
find its spectrum, eigenfunctions, inversion, invariants 
(invariant $q$-difference operators), and give its interpretation 
as a fractional $q$-integration operator. We also give expansions
of the $A_1$ Macdonald polynomials into the eigenfunctions of the 
separating operator $M$ and vice versa.
\vskip 1cm
\newpage
\newsection{Introduction}
One of the most powerful methods for solving the spectral problem
of a quantum integrable system is that of {\it Separation of Variables}
(SoV). Given a quantum integrable
system with $n$ degrees of freedom defined by a complete set of 
$n$ commuting Hermitian operators (integrals of motion) $H_j$
acting in some Hilbert space ${\cal H}$,
\beq
 [H_j,H_k]=0, \qquad j,k=1,2,\ldots,n\,,
\label{eq:comm-H}
\eeq
one tries to find a linear operator $M$ sending any common eigenvector 
$P_{\l}$ of the $H_j$'s,
\beq
 H_j\;P_{\l}=h_{\l;j} \;P_{\l}\,,
\label{eq:spec-H}
\eeq
labelled by the quantum numbers
$\l=\{\l_1,\ldots,\l_n\}$, into the product $F_\l$
\beq
 M: P_{\l}\rightarrow F_\l=\prod_{j=1}^n f_{\l;j}(y_j)
\label{eq:def-SoV}
\eeq
of functions $f_{\l;j}(y_j)$ of one variable each
referred to as partial, or separated, or factorized
functions. Notice that the operator
$M$ does not depend on the spectrum 
$\{h_{\l;j}\}_{j=1}^n$ (or the quantum numbers $\l$)
 and is thought of as an intertwiner 
between two representations of the eigenfunctions of the integrals
$H_j$, namely: 
the initial representation and the separated (or factorized)
$y$-representation. The essence of the method is that the original
{\it multivariable} spectral problem (\ref{eq:spec-H})
is transformed by such a transformation 
into a set of (simpler) {\it multiparameter} spectral problems, each being in 
one variable, given by the separation equations of the form
\beq
 {\cal D}_j\left(y_j,\frac{\partial}{\partial y_j};
h_{\l;1},\ldots,h_{\l;n}\right)f_{\l;j}(y_j)=0\,,\qquad j=1,\ldots,n\,,
\label{eq:sep-eq-gen}
\eeq
where ${\cal D}_j$ are usually some differential
or finite-difference operators in variable $y_j$ depending on the
spectral parameters $h_{\l;k}$. In the context of the classical Hamiltonian
mechanics the above construction corresponds precisely to the 
standard definition of SoV in the Hamilton-Jacobi equation.

For a recent review of the SoV method see \cite{S1}.
See also \cite{KS1,KS2} as for applications of SoV to the 
$3$-particle Calogero-Sutherland model and its
relativistic (or $q$-) analog, respectively. 
We would like to stress here that our definition of SoV is
not restricted by purely coordinate changes of variables,
as it is usually (historically) supposed 
(see, for instance, the book \cite{Ka} and references
therein). Instead, we allow the separating transform to be
a generic canonical transformation in classical mechanics
or a linear (integral) operator in quantum case.

The general multiparameter spectral theory for systems
of equations like (\ref{eq:sep-eq-gen}) including, in particular,
questions related to completeness of eigenfunctions, the Parseval
equality and the like, has been studied extensively
\cite{ARS1,Atk1,Atk2,SLE,SLE1,V}. 
The book \cite{Atk2} (see also
survey article \cite{Atk1}) gives a comprehensive treatment of 
multiparameter spectral theory in finite dimensional cases. 
Various results concerning unbounded (e.g. differential) operators
are brought together in the
monograph \cite{SLE} where the author provides a study of
the infinite dimensional case via the theory of several
commuting operators in Hilbert space [$H_j^{(y)}$]. 

In the present paper we will not concentrate on the multiparameter
spectral problem because for the particular system we study
we do know explicitly the spectrum and the solution of the separation
equation. Although, when doing SoV for more complicated integrable
systems, like Toda lattice or elliptic Calogero-Moser problem, we expect 
the above-mentioned  theory will play an important role. 

The main problem of SoV (apart from the study of the separation
equation and its solution) is to try to describe the separating 
operator $M$ (and its inversion) in the most explicit terms. 
The eigenvectors $P_\l$ in a particular representation
become special functions or orthogonal polynomials 
in many variables (for instance, Jack or Macdonald polynomials)
while $f_{\l;j}$ in (\ref{eq:def-SoV}) become
special functions in one variable. 
In this interpretation the corresponding kernel (nucleus)
$\M$ of the integral 
operator $M$ is some special function in $2n$ variables and the method of SoV 
provides us with two fundamental integral {\it relations} between special
functions in one and many variables:
\begin{subeqnarray}
\prod_{i=1}^nf_{\l;i}(y_i)&=&\int\cdots\int dx_1\ldots dx_n \;\M(y|x)\;
P_\l(x_1,\ldots,x_n) \\
P_\l(x_1,\ldots,x_n)&=&\int\cdots\int dy_1\ldots dy_n 
\;\tilde\M(x|y)\;\prod_{i=1}^n f_{\l;i}(y_i)\,.
\label{eq:2-relations}
\end{subeqnarray}
Here $\tilde\M(x|y)$ represents the kernel of the inverse $M^{-1}$
and, we repeat, both kernels do not depend on $\l$. 

The first relation in (\ref{eq:2-relations}) looks like a ``product formula''
for the one-variable functions $f_{\l;j}(y_j)$ while the second one
gives an integral representation 
for the special function $P_\l(x_1,\ldots,x_n)$
of many variables in terms of the special functions, $f_{\l;j}(y_j)$, of one
variable. These two types of integral relations are not entirely new 
in the theory of special functions. In particular, when $n=2$ both integral
relations have appeared before in one-variable theory. 
Below we give a brief review of the literature on integral 
relations/equations (\ref{eq:2-relations}) for the special functions and 
orthogonal polynomials in one variable. Let us start with the relations of
2-nd type (\ref{eq:2-relations}b).

In 1846 Liouville \cite{L} demonstrated that Lam\'{e} polynomials
satisfy certain non-linear integral equations of the form
\beq
f(x)=\int\int dy_1 \,dy_2 \;\tilde\M(x|y)\;f(y_1)\;f(y_2)\,.
\label{xxxx}\eeq
However, he failed to specify for which of the eight types 
of Lam\'{e} polynomials these equations were valid. Later on  
Sleeman \cite{SLE2} showed, using the results obtained by Arscott
in \cite{ARS2} (see also the work \cite{ARS1} and the book \cite{ARS}), 
that integral equations and 
relations may be constructed which are satisfied by all eight types
of Lam\'{e} polynomials and Lam\'{e} functions of the second type. 
The kernels (nuclei) appearing in all these equations are the ``potential''
Green's functions. These results were further generalized in \cite{SLE2}
to construction of integral equations and relations satisfied
by ellipsoidal wave functions of the first and third kinds, the kernels
now being the ``free space'' Green's functions. The partial differential
equations for the kernels $\tilde\M(x|y)$ as well as the corresponding
equality of the form (\ref{xxxx}) for these cases
were first written down in \cite{ARS1} (see also the book \cite{ARS}) 
under the names ``Second integral theorem'' and ``Second integral 
equation''. Actually, the very existence of the relation (\ref{xxxx})
is closely connected to the corresponding 2-parameter spectral 
problem for the Lam\'{e} polynomials/functions and for the ellipsoidal
wave functions as was first realized by Arscott \cite{ARS}. This idea 
was generalized further to the case of $n$-parameter spectral problem
where one gets the following equation 
\beq
f(x)=\int\cdots\int dy_1\,\cdots\,dy_n \;\tilde\M(x|y)\; \prod_{j=1}^nf(y_j)
\label{5x}\eeq
as well as equation of the form
\beq
\prod_{i=1}^nf(x_i)
=\int\cdots\int dy_1\,\cdots\,dy_n \;\tilde\M(x|y)\; \prod_{j=1}^nf(y_j)\,.
\label{x5}\eeq
An abstract formulation (in Hilbert space) of both above-mentioned
equations
(and also the theory of solvability of some operator equations behind
it) and the characteristic (determinantal) equations for the 
corresponding kernels $\tilde\M(x|y)$ can be found in the monograph
\cite{SLE} (Chapter 6 titled ``An abstract relation''), 
see also the works \cite{ks,Br}. For the case when $f(y)$
is the ellipsoidal wave function the relations (\ref{5x}) and (\ref{x5}) 
(with $n=2$) were derived for the first time by Malurkar \cite{Mal}
and M\"oglich \cite{Mog}, respectively. 

Let us proceed now to the 1-st type integral relations (\ref{eq:2-relations}a).
The particular form of these relations
\beq
f(y_1)\,f(y_2)=\int dx \;\M(y|x)\;f(x)
\label{6x}\eeq
is known in the literature
as product formulas for the orthogonal polynomials 
and special functions in one variable. 
The first example was provided by Gegenbauer \cite{Ge}
for his polynomials. This approach was developed further for several
other classical special functions and orthogonal polynomials
\cite{Vi,Sha,Ga,Ko,AbK,RV,Ra,GR,CS,CMS1,Ma} and  
culminated in \cite{CMS2} where the product formulas for
angular and radial spheroidal wave functions were presented.
These formulas 
generalize the results for the Gegenbauer (ultraspherical)
polynomials and functions as well as for the Mathieu functions. 
The modern approach to the construction of some of the corresponding 
kernels $\M(y|x)$ employs a PDE technique based on Riemann's
integration method to solve a Cauchy problem for an associated 
hyperbolic differential equation (cf. \cite{CMS1,Ma,CMS2}). 

Product formulas are very useful identities associated with a discrete
or continuous function system. Generally, such a formula represents 
the product of any two values of an orthogonal polynomial in terms 
of a Stieltjes integral which depends linearly on the polynomial
itself. It is a very important tool in the harmonic analysis of orthogonal
expansions where it usually defines a convolution product on a function
space (this space 
usually becomes a Banach algebra under this operation
and the convolution product 
plays the same role as ordinary convolution in Fourier 
analysis). Notice that the right hand side
of a product formula usually defines
a generalized translation operator and that this
in turn can be used to 
establish an appropriate convolution structure on the space and often
a hypergroup structure on the corresponding interval (cf. \cite{CMS1}
for details). The derivation of new product formulas yields new
convolution structures and hypergroups; moreover it may be used to 
obtain more information about the special functions $f$ themselves,
since they are usually not given in closed form as soon as they
lie in the ``land beyond Bessel'' \cite{ARS3}. 

In 1914 Whittaker \cite{W} conjectured that the Heun 
differential equation \cite{BE,He} is the simplest equation of Fuchsian type
whose solution cannot be represented by a contour integral
(of a simpler, elementary function); instead the nearest 
approach to such a solution is to find an integral equation 
satisfied by a solution of the differential equation. 
This is indeed a general situation for the higher special
functions (cf. \cite{ARS3}) where one could not hope to
get as many explicit formulas and representations as there are
for special functions of the hypergeometric type. For instance,
the technique of series solutions leads 
to a higher order recursion relation involving at least
three successive coefficients (for functions of the Heun type).
This observation applies not only to power series but to series 
of other forms. Although there is no proof that a two-term recursion
cannot be obtained. Usually, no integral expressions for higher special 
functions are available; instead they satisfy certain integral equations
like the ones mentioned above. 

In the present paper we aim to demonstrate that the integral equations 
similar to those
just discussed appear very naturally in the unified approach 
of SoV. One (among the many) way to introduce a multivariable special
function generalizing a given one-variable special function is to associate
it to a quantum integrable family, i.e. to obtain it as a common eigenfunction
of a commutative ring of $n$ quantum integrals of motion. In such a manner
many known functions can be described and characterized. 
Assuming the hypothesis
that any quantum integrable system can be separated we arrive at
an interesting problem of studying various integral operators, in an 
attempt to factorize a particular function or polynomial. That was our leading 
motivation for studying $3$-variable Jack polynomials in \cite{KS1} 
and $3$-variable Macdonald polynomials in \cite{KS2} which are symmetric
orthogonal polynomials associated to the root system $A_2$ characterized
by being common eigenfunctions of the corresponding quantum 
Calogero-Sutherland \cite{Ca,Su} and Ruijsenaars \cite{Ru} models
of the trigonometric type.
It appears that there are explicit separating operators $M$
for both problems which factorize polynomials and which can be explicitly 
inverted giving new integral representations for the related polynomials.

Here we apply the SoV technique
to even simpler $2$-variable case of the $A_1$ Macdonald polynomials. 
In particular, we give a new proof of the ``quadratic'' integral 
equations for the continuous $q$-ultraspherical
polynomials. The separating integral operator $M$ expressed in terms
of the Askey-Wilson operator is studied in detail: apart from writing
down the characteristic (``separation'') equations it satisfies, we
find its spectrum, eigenfunctions, inversion, invariants 
(invariant $q$-difference operators), and give its interpretation 
as a fractional $q$-integration operator. We also give expansions
of the $A_1$ Macdonald polynomials into the eigenfunctions of the 
separating operator $M$ and vice versa. The structure of the operators $M$
for $A_1$ and $A_2$ cases turns out to be quite similar which indicates an 
intimate relation between the two SoV problems.

The structure of the paper is as follows. In Section 2 we collect
basic information about continuous $q$-ultraspherical polynomials 
including the product formula in integral form. In Section 3
we give the main properties of the Askey-Wilson type integral operator
$M_{\a\b}$ which provides us with an important tool for construction 
and description of the separating operator $M$ in the key Section, 4.
In the next Section, 5, we derive expansions for the $A_1$ 
Macdonald polynomials in terms of the eigenfunctions of the 
separating operator $M$, as well as the dual expansions.
As a bi-product we obtain an interesting interpretation of $M$ as a
quantum integrable map (discrete-time system) and find its invariants
(integrals of motion): the $q$-difference operators  $N_1$ and $N_2$. 
We provide some concluding
remarks in the final Section, 6, and give an Appendix
where the classical SoV for the trigonometric 
$A_1$ type Ruijsenaars model is performed
and a crucial (for the present paper) inter-relation 
between the $A_2$ and $A_1$ problems is explained in detail. 

\newsection{Continuous $q$-ultraspherical polynomials}
The main references for the formulas and notations 
in this Section are \cite{GR,Koe}. 

Continuous $q$-ultraspherical polynomials $C_n(\xi;\b|q)$ \cite{AI} 
is an important class of orthogonal polynomials of basic
hypergeometric type within the Askey scheme. They
can be obtained from the general Askey-Wilson polynomials \cite{AW} 
through a specification of their four parameters, so that $C_n(\xi;\b|q)$
depend only on one parameter $\b$,  apart from the degree $n$ and
basic parameter $q$. They are defined through the (terminated)
${}_2\oldphi{}_1$ basic hypergeometric series
\bea
C_n(\xi;\b|q)&=&\frac{(\b;q)_n}{(q;q)_n}\;\e^{in\theta}\;
{}_2\oldphi{}_1\left[\matrix{q^{-n},\,\b\cr
\b^{-1}q^{1-n}};q,q\b^{-1}\e^{-2i\theta}\right]\nonumber\\
&=&\sum_{k=0}^n\frac{(\b;q)_k(\b;q)_{n-k}}{(q;q)_k(q;q)_{n-k}}
\;\e^{i(n-2k)\theta}\,,\qquad \xi=\cos\theta\,.
\label{2.1}\eea
Notice that in the limit $q\uparrow 1$ we get usual Gegenbauer
(or ultraspherical) polynomials $C_n^\l(\xi)$ (see \cite{BE,Koe})
\beq
\lim_{q\uparrow 1}C_n(\xi;q^\l|q)=\sum_{k=0}^n
\frac{(\l)_k(\l)_{n-k}}{k!(n-k)!}
\;\e^{i(n-2k)\theta}=C_n^\l(\xi)\,.
\label{2.2}\eeq
Polynomials $C_n(\xi;\b|q)$ are orthogonal polynomials on $[-1,1]$
with the following measure ($|q|<1, \; |\b|<1$):
\bea
&&\frac{1}{2\pi}\int_{-1}^1C_m(\xi;\b|q)\;C_n(\xi;\b|q)\;
\left|\frac{(\e^{2i\theta};q)_\infty}{(\b\e^{2i\theta};q)_\infty}
\right|^2\frac{d\xi}{\sqrt{1-\xi^2}}\nonumber\\
&&\qquad =
\frac{(\b,\b q;q)_\infty}{(\b^2,q;q)_\infty}
\frac{(\b^2;q)_n}{(q;q)_n}\frac{(1-\b)}{(1-\b q^n)}\;\delta_{mn}\,.
\label{2.3}
\eea
They satisfy the recurrence relation
\beq
2(1-\b q^n)\,\xi\, C_n(\xi;\b|q)=
(1-q^{n+1})\,C_{n+1}(\xi;\b|q)+(1-\b^2q^{n-1})\,C_{n-1}(\xi;\b|q)
\label{2.4}\eeq
and $q$-difference equation of the form
\beq
(1-q)^2D_q\,[w(\xi;q^\half\b|q)\,D_qy(\xi)]+\l_n\,w(\xi;\b|q)\,y(\xi)=0\,,\quad
y(\xi)=C_n(\xi;\b|q)\,,
\label{2.5}\eeq
where 
\beq
w(\xi;\b|q):=\left|\frac{(\e^{2i\theta};q)_\infty}
{(\b\e^{2i\theta};q)_\infty}\right|^2\,,\qquad 
\l_n=4q^{-n+1}(1-q^n)(1-\b^2q^n)
\label{2.6}\eeq
and $D_qf(\xi):=\frac{\delta_qf(\xi)}{\delta_q \xi}$ with $\delta_q
f(\e^{i\theta})=f(q^\half\e^{i\theta})-f(q^{-\half}\e^{i\theta})$,
$\; \xi=\cos\theta$. 

There is a simple generating function for these polynomials:
\beq
\frac{(\b\e^{i\theta}z,\b\e^{-i\theta}z;q)_\infty}
{(\e^{i\theta}z,\e^{-i\theta}z;q)_\infty}=
\sum_{n=0}^\infty C_n(\xi;\b|q)\;z^n\,, \qquad \xi=\cos\theta\,.
\label{2.7}\eeq
Finally, let us mention the product formula in integral form 
(see (8.4.1)--(8.4.2) in \cite{GR} and work \cite{RV})
\beq
C_n(\xi;\b|q)\;C_n(\eta;\b|q)=\frac{(\b^2;q)_n}{(q;q)_n}\;\b^{-n/2}
\int_{-1}^1K(\xi,\eta,\zeta;\b|q)\;C_n(\zeta;\b|q)\;d\zeta
\label{2.8}
\eeq
where ($\xi=\cos\theta$, $\eta=\cos\oldphi$, $\zeta=\cos\psi$)
\bea
K(\xi,\eta,\zeta;\b|q)&=&
\frac{(q,\b,\b;q)_\infty\abs{(\b\e^{2i\theta},\b\e^{2i\oldphi};q)_\infty}^2}%
{2\pi(\b^2;q)_\infty\sqrt{1-\zeta^2}}\,\nonumber\\
&\times& w(\e^{i\psi};\b^\half\e^{i\theta+i\oldphi},
\b^\half\e^{-i\theta-i\oldphi},\b^\half\e^{i\theta-i\oldphi},
\b^\half\e^{i\oldphi-i\theta})\,.
\label{2.9}
\eea
Function $w(x;a,b,c,d)$ in (\ref{2.9}) is defined as follows:
\beq
 w(x;a,b,c,d):=
\frac{(x^2,x^{-2};q)_\infty}%
{(ax,ax^{-1},bx,bx^{-1},cx,cx^{-1},dx,dx^{-1};q)_\infty}\,.
\label{eq:AW-kernel}
\eeq
In the Section 4 we will give a new proof of this product formula. 
\vskip 0.2cm

\noindent
\begin{rem}
It is always possible to give a general expression for the kernel $K$
as an infinite series. Indeed,
consider a product formula in integral form for some polynomials
$p_n(x)$: 
\beq
p_n(x)\;p_n(y)=\int_a^bdz\;K(x,y,z)\;p_n(z)
\label{ssss}\eeq
and let the orthogonality relation for the polynomials be given by
\beq
\int_a^b p_n(x)\;p_m(x)\;w(x)\;dx=b_n\;\delta_{nm}\,.
\label{ffff}\eeq
Then, the left hand side $p_n(x)\,p_n(y)$ of (\ref{ssss}) 
can be thought of as the expression for the 
Fourier coeficients  in the expansion of the function
$K(x,y,z)$ into the basis of orthogonal polynomials $p_n(z)$. 
Hence, we always have the formula for the kernel $K(x,y,z)$ in terms
of the product of 3 polynomials $p_n(x)\,p_n(y)\,p_n(z)$ as a consequence
of (\ref{ssss}) and (\ref{ffff}):
$$
K(x,y,z)=w(z)\;\sum_{n=0}^\infty\frac{p_n(x)\,p_n(y)\,p_n(z)}{b_n}\,.
$$
Applying this to the kernel \Ref{2.9} and using \Ref{2.3} we obtain
the expansion
$$
K(\xi,\eta,\zeta;\b|q)=\frac{(q,\b^2;q)_\infty}
{2\pi (\b,\b;q)_\infty\sqrt{1-\zeta^2}}\;
\left|\frac{(\e^{2i\psi};q)_\infty}
{(\b\e^{2i\psi};q)_\infty}\right|^2\;\times
$$
\beq
\times\;\sum_{n=0}^\infty
\;\b^{n/2}\;(1-\b q^n)\;\left(\frac{(q;q)_n}{(\b^2;q)_n}\right)^2\;
C_n(\xi;\b|q)\,C_n(\eta;\b|q)\,C_n(\zeta;\b|q)\,.
\label{dddd}\eeq
\label{rem:cubic-series}
\end{rem}
\vskip 0.2cm

\noindent
\begin{rem}
It was shown in \cite{RV} that in the limit $q\uparrow 1$ the product
formula (\ref{2.8}) goes into the product formula due to Gegenbauer
\cite{Ge} 
\beq
\frac{C_n^\l(x)\;C_n^\l(y)}{C_n^\l(1)\;C_n^\l(1)}
=\int_{-1}^1K(x,y,z)\;\frac{C_n^\l(z)}{C_n^\l(1)}\;dz\,,
\label{7.1}\eeq
where
\be
K(x,y,z)=\left\{\begin{array}{rcl}
\dfrac{(1-x^2-y^2-z^2+2xyz)^{\l-1}}%
{B(\l,{\scriptstyle\half})[(1-x^2)(1-y^2)]^{\l-\half}}
\,,&\quad& 1-x^2-y^2-z^2+2xyz>0\,,\cr
0\,,&\quad& 1-x^2-y^2-z^2+2xyz\leq 0\,.\end{array}
\right.\ee
\end{rem}
\vskip 0.2cm

Now we describe the relation between the $q$-ultraspherical polynomials and 
$A_1$ Macdonald polynomials.
Let $\Symx$ be the space of symmetric Laurent polynomials in two
variables $x_1$, $x_2$.
Let us introduce the polynomials $P_\l\in\Symx$, 
\beq
 P_\l(x_1,x_2)=x_+^{\abs{\l}}\;
\frac{(q;q)_{\l_{21}}}{(t;q)_{\l_{21}}}\;
C_{\l_{21}}(\tfrac{x_-+x_-^{-1}}{2}\,;t | q)\,, \qquad
 x_\pm\equiv(x_1x_2^{\pm1})^{\half}\,,
\label{eq:def-Pl}
\eeq
depending on real parameters $q,t\in(0,1)$ 
and indexed by the pairs of integers $\l\equiv\{\l_1,\l_2\}\in\Z^2$ such that
$\l_1\leq\l_2$. We use the notation
 $\abs{\l}\equiv\l_1+\l_2$,
$\l_{21}\equiv\l_2-\l_1$, $t=q^g$, $g>0$.

\begin{prop}
The polynomials $P_\l(x_1,x_2)$ (\ref{eq:def-Pl})
coincide with the standard Macdonald polynomials for the root system $A_1$.
\end{prop}
\Proof
{}From \Ref{2.5}--\Ref{2.6} it follows that $P_\l$ are the eigenfunctions
\beq
 H_j\;P_\l=h_{\l;j}\;P_\l
\eeq
of two finite-difference operators in the space ${\rm Sym}(x_1,x_2)$	
\beq
 H_1=v_{12}\;T_{q,x_1}+v_{21}\;T_{q,x_2}\,, \qquad H_2=T_{q,x_1}\,T_{q,x_2}
\label{eq:def-H}
\eeq
where
\beq
 v_{jk}=\frac{t^{\half}x_j-t^{-\half}x_k}{x_j-x_k}
\label{eq:def-w}
\eeq
and $T_{q,x}$ is the $q$-shift operator in the variable $x$
\beq
 T_{q,x}f(x)=f(qx)\,.
\eeq
The corresponding eigenvalues $h_{\l;j}$ are
\beq
 h_{\l;1}=t^{-\half}q^{\l_1}+t^{\half}q^{\l_2}\,, \qquad
 h_{\l;2}=q^{\abs\l}\,.
\eeq
Introducing the basis of monomial symmetric
functions $m_\l$ in $\Symx$
\beq
 m_\l(x_1,x_2)=\left\{\begin{array}{rcl}
      x_1^{\l_1}x_2^{\l_2}+x_1^{\l_2}x_2^{\l_1}\,,&\quad& \l_1<\l_2 \\
      x_1^{\l_1}x_2^{\l_2}\,, &\quad& \l_1=\l_2
 \end{array}\right.
\label{eq:def-monomial}
\eeq
and using the definition \Ref{2.1} we observe that $P_\l$ expands in 
terms of $m_\nu$ as
\beq
 P_\l=\sum_{\scriptstyle \nu\prec\l \atop \scriptstyle \abs{\nu}=\abs{\l}}
u_{\l\nu}\;m_\nu\,,
\eeq
\beq
u_{\l\nu}=\frac{(q;q)_{\l_{21}}}{(t;q)_{\l_{21}}}
\frac{(t;q)_{\nu_1-\l_1}}{(q;q)_{\nu_1-\l_1}}
\frac{(t;q)_{\l_2-\nu_1}}{(q;q)_{\l_2-\nu_1}}\,.
\label{expan}\eeq
The ordering $\prec$ is defined here as
\beq
   \nu\prec\l \Longleftrightarrow \l_1\leq\nu_1\leq\nu_2\leq\l_2\,.
\label{eq:order}
\eeq
{}From \Ref{2.1} there follows also the normalization condition 
$u_{\l\l}=1$.
The enlisted properties of $P_\l$ constitute precisely the definition of
$A_1$ Macdonald polynomials (see \cite{Mac}, Chapter VI).
\endproof

Let us introduce also the separated polynomials $f_\l(y)$
\bea
 f_\l(y)&=&y^{\l_1}\, {}_2\oldphi_1\left[
\matrix{t,q^{-\l_{21}} \cr t^{-1}q^{1-\l_{21}} };q,t^{-2}qy\right] 
\nonumber \\
&=&y^{\frac{\abs{\l}}{2}}t^{-\frac{\l_{21}}{2}}
\;\frac{(q;q)_{\l_{21}}}{(t;q)_{\l_{21}}}\;
C_{\l_{21}}(\cos\theta;t | q)\,, \qquad y=t\,e^{-2i\theta}\,.
\label{eq:def-fl}
\eea
It follows from \Ref{2.5} and \Ref{2.1} that $f_\l(y)$ satisfies the 
$q$-difference equation
\beq
 [\;t(1-qy)\,T_{q^2,y}-t^{\half}(t-qy)\,h_{\l;1}\,T_{q,y}
+(t^2-qy)\,h_{\l;2}\;]\;f_\l(y)=0
\label{eq:separ-eq}
\eeq
and expands in $y$ as
\beq
 f_\l(y)=\sum_{k=\l_1}^{\l_2}\chi_{k}\;y^k, \qquad
\chi_{k}=(t^{-2}q)^{k-\l_1}
\frac{(t,q^{-\l_{21}};q)_{k-\l_1}}{(q,t^{-1}q^{1-\l_{21}};q)_{k-\l_1}}\,.
\eeq
In particular,
\beq
 \chi_{\l_1}=1\,, \qquad \chi_{\l_2}=t^{-\l_{21}}\,.
\label{eq:chi-l12}
\eeq
\vskip 0.2cm

\noindent
\begin{rem} In \cite{KS2} a general formula is given for the separated 
polynomials for any root system $A_{n-1}$. Putting $n=2$ in the formula (5.1)
from \cite{KS2} one obtains another expression for $f_\l(y)$:
\beq
 f_\l(y)=y^{\l_1}(y;q)_{1-2g}\;{}_2\oldphi_1
\left[\matrix{t^{-2}q^{1-\l_{21}},t^{-1}q \cr 
        t^{-1}q^{1-\l_{21}}};q,y\right].
\label{eq:fl-alt-expr}
\eeq
The equivalence of \Ref{eq:def-fl} and \Ref{eq:fl-alt-expr}
follows from the $q$-analog of
Euler's transformation  formula \cite{GR}, (1.4.3).
\label{rem:Euler-tr}
\end{rem}
\vskip 0.2cm

The formula \Ref{2.8} can now be interpreted as the one expressing SoV  for
the pair of integrals of motion \Ref{eq:def-H} for the quantum integrable
system with 2 degrees of freedom (which is the $A_1$ type trigonometric
Ruijsenaars' system \cite{Ru}).
\vskip 0.2cm

\noindent
\begin{rem}  As one can see, the $A_1$ Macdonald polynomials $P_\l(x_1,x_2)$
\Ref{eq:def-Pl} are already in the factorized form when written in terms 
of the variables $x_\pm$. This corresponds to another (trivial) 
SoV which is purely coordinate one (local transform):
$(x_1,x_2)\mapsto (x_+,x_-)$, as opposed to the integral separating 
transform \Ref{2.8}. This is a simple demonstration of the fact that SoV
for a given integrable system (a given multi-variable special function)
might be not unique in general. 
\label{rem:center-of-mass}
\end{rem}

\newsection{Operator $M_{\a\b}$ and its properties}
In this and the next sections we introduce and study 
the integral operator $M_\xi$ performing 
the separation of variables in the $A_1$ Macdonald polynomials \Ref{eq:def-Pl}.
The operator is  most conveniently described in terms
of a slightly more general integral operator $M_{\a\b}$ which, in turn,
is closely related to the fractional $q$-integration operator $I^\a$.

Our main technical tool is the famous Askey-Wilson integral identity 
\cite{AW,GR}
\beq
 \frac{1}{2\pi i}\int\limits_{\G_{abcd}}\frac{dx}{x}\,w(x;a,b,c,d)=
\frac{2(abcd;q)_\infty}{(q,ab,ac,ad,bc,bd,cd;q)_\infty}
\label{eq:AW-integral}
\eeq
where $w$ is given by \Ref{eq:AW-kernel}.
The cycle $\G_{abcd}$ depends on complex parameters $a$, $b$, $c$, $d$ and is 
defined as follows. Let $C_{z,r}$ be the counter-clockwise oriented circle 
with the center $z$ and radius $r$. 
If $\mbox{max}(\abs{a},\abs{b},\abs{c},\abs{d},\abs{q})<1$ 
then $\G_{abcd}=C_{0,1}$.
The identity \Ref{eq:AW-integral}
can be continued analytically for the values of parameters $a$, $b$, $c$,
$d$ outside the unit circle provided the cycle $\G_{abcd}$ is deformed 
appropriately. In general case 
\beq
 \G_{abcd}=C_{0,1}+\sum_{z=a,b,c,d}
\sum_{\scriptstyle k\geq0\atop \scriptstyle \abs{z}q^k>1}
(C_{zq^k,\eps}-C_{z^{-1}q^{-k},\eps})\,,
\eeq
$\eps$ being small enough for $C_{z^{\pm1}q^{\pm k},\eps}$ to encircle 
only one pole of the denominator. 

Now put
\beq
 a=yq^{\frac\a2}\,, \qquad b=y^{-1}q^{\frac\a2}\,, \qquad
 c=rq^{\frac\b2}\,, \qquad d=r^{-1}q^{\frac\b2}\,, 
\label{eq:def-abcd}\eeq
and introduce the notation 
\beq
 \Lq\nu xy:=(\nu xy,\nu xy^{-1},\nu x^{-1}y,\nu x^{-1}y^{-1};q)_\infty\,.
\eeq
The kernel
\beq
 \M_{\a\b}(r,y| x)=
\frac{(1-q)(q;q)_\infty^2\,(x^2,x^{-2};q)_\infty
\;\Lq{q^{\frac{\a+\b}{2}}}ry }%
{2B_q(\a,\b)
\;\Lq{q^{\frac\a2}}yx\,\Lq{q^{\frac\b2}}rx  }
\label{eq:ker-Mab}
\eeq
defines the integral operator 
\beq
 (M_{\a\b}^rf)(y):=\frac{1}{2\pi i}\int\limits_{\Gamma_{\a\b}^{ry}}
\frac{dx}{x}\;\M_{\a\b}(r,y|x)\;f(x)\,,
\label{eq:def-Mab}
\eeq
the contour $\G_{\a\b}^{ry}$ being obtained from $\G_{abcd}$ 
by the substitutions (\ref{eq:def-abcd}). Notice that in 
\Ref{eq:ker-Mab}--\Ref{eq:def-Mab} we consider $x$ and $y$ as
arguments and $r$ as a parameter. 

The following properties of the operator $M_{\a\b}^r$ are proved in
\cite{KS2}, Appendix B.

\begin{prop}
\hfill\par
\begin{enumerate}
\item 
Let ${\rm Refl}(x)$ be the space of reflexive 
(invariant w.r.t.\ $x\rightarrow x^{-1}$) Laurent polynomials in $x$. 
Then $M_{\a\b}^r:{\rm Refl}(x)\rightarrow{\rm Refl}(y)$.

\item 
In particular, from the Askey-Wilson identity \Ref{eq:AW-integral}
it follows immediately that $M:1\rightarrow1$.

\item 
More generally, consider a Laurent polynomial 
$R^{\a\b}_{j_1j_2k_1k_2}(r,x)\in{\rm Refl}(x)$, where 
$j_{1,2},k_{1,2}\in\Z_{\geq0}$
\begin{eqnarray}
 R^{\a\b}_{j_1j_2k_1k_2}(r,x)&:=&
   (q^{\frac{\a}{2}}yx,q^{\frac{\a}{2}}yx^{-1};q)_{j_1}\;
(q^{\frac{\a}{2}}y^{-1}x,q^{\frac{\a}{2}}y^{-1}x^{-1};q)_{j_2}
\nonumber \\
 &\times&(q^{\frac{\b}{2}}rx,q^{\frac{\b}{2}}rx^{-1};q)_{k_1}\;
(q^{\frac{\b}{2}}r^{-1}x,q^{\frac{\b}{2}}r^{-1}x^{-1};q)_{k_2}\,.
\label{eq:def-R}
\end{eqnarray}
Then
\begin{eqnarray}
M_{\a\b}^r: R^{\a\b}_{j_1j_2k_1k_2}(r,x)&\rightarrow&
\frac{(q^\a;q)_{j_1+j_2}\,(q^\b;q)_{k_1+k_2}}%
{(q^{\a+\b};q)_{j_1+j_2+k_1+k_2}} 
\nonumber \\
&&\times
(q^{\frac{\a+\b}{2}}ry;q)_{j_1+k_1}\;
(q^{\frac{\a+\b}{2}}ry^{-1};q)_{j_2+k_1} \nonumber \\
&&\times
(q^{\frac{\a+\b}{2}}r^{-1}y;q)_{j_1+k_2}\;
(q^{\frac{\a+\b}{2}}r^{-1}y^{-1};q)_{j_2+k_2}\,.
\label{eq:actMR}
\end{eqnarray}

The last formula allows to calculate effectively the action of
$M_{\a\b}^r$ on any polynomial from ${\rm Refl}(x)$.

\item The inversion of $M_{\a\b}^r$ is given by the formula
\beq
 \left(M_{\a\b}^r\right)^{-1}=M_{-\a,\a+\b}^r\,,
\label{3.9}\eeq
the corresponding kernel being
\beq
 \tilde\M_{\a\b}^r(x| y)=
\frac{(1-q)(q;q)_\infty^2\;(y^2,y^{-2};q)_\infty
\;\Lq{q^{\frac{\b}{2}}}rx  }%
{2B_q(-\a,\a+\b)
\;\Lq{q^{-\frac{\a}{2}}}yx\;\Lq{q^{\frac{\a+\b}{2}}}ry   }\,.
\label{3.10}\eeq

\item 
The operator $M_{\a\b}^r$ simplifies drastically when either of the
parameters $\a$, $\b$ takes negative integer values.
Let $\a=-g$, $g\in\Z_{\geq0}$. Then $M_{\a\b}^r$ turns into the 
$q$-difference operator of order $g$:
\beq
 M_{-g,\b}^r: f(x)\rightarrow \sum_{k=0}^{g} \xi_k(r,y)\;f(q^{k-\frac{g}{2}}y)
\label{eq:M_(-gb)}
\eeq
where
\begin{eqnarray}
\lefteqn{
 \xi_k(r,y)=(-1)^k q^{-\frac{k(k-1)}{2}} 
\left[\begin{array}{c}g \\ k\end{array}\right]_q y^{-2k}\;
(1-q^{g-2k}y^{-2}) }\nonumber \\
&&\times
\frac{(q^{\frac{\b-g}{2}}ry,q^{\frac{\b-g}{2}}r^{-1}y;q)_k\;
(q^{\frac{\b-g}{2}}ry^{-1},q^{\frac{\b-g}{2}}r^{-1}y^{-1};q)_{g-k}}%
{(q^{\b-g};q)_{g}\;(q^{-k}y^{-2};q)_{g+1} }\,.
\label{def-xi-k}\end{eqnarray}
\end{enumerate}
\label{th:properties-of-Mab}
\end{prop}

It was mentioned briefly in \cite{KS2} 
that, similarly to the nonrelativistic case \cite{KS1}, 
the operator $M_{\a\b}^r$ has a natural 
interpretation in terms of a {\it fractional $q$-integration} operator.
To clarify this remark,
take the formula \Ref{eq:ker-Mab} for the kernel 
$\M_{\a\b}(r,y| x)$
and rewrite it in the form
\beq
 \M_{\a\b}(r,y| x)=\frac{\psi_r^{\b-1}(x)}{\psi_r^{\a+\b-1}(y)}
\;\I^\a(r,y| x)\,,
\eeq
with
\beq
 \psi_r^\nu(x):=
\frac{\Lq{q^{\half}}{r}{x}}{\Gamma_q(\nu+1)\;\Lq{q^{\frac{\nu+1}{2}}}{r}{x}}\,,
\label{psi}\eeq
where the kernel $\I^\a(r,y| x)$ is defined as follows:
\beq
 \I^\a(r,y| x)=
\frac{(1-q)(q;q)_\infty^2\;(x^2,x^{-2};q)_\infty\;\Lq{q^{\half}}{r}{y}}%
{2\Gamma_q(\a)\;\Lq{q^{\frac{\a}{2}}}{y}{x}\;\Lq{q^{\half}}{r}{x}}\,.
\label{qu-qu}\eeq
In the operator form we have
\beq
 M_{\a\b}^r=\frac{1}{\psi_r^{\a+\b-1}}\circ I^\a\circ\psi_r^{\b-1}
\label{eq:M-I}
\eeq
where $\psi_r^\nu$ are thought of as multiplication operators
and the operator $I^\a$ corresponds to the kernel
$\I^\a(r,y| x)$ \Ref{qu-qu}
\beq
 I^\a:f(x)\rightarrow\frac{1}{2\pi i}\int\limits_{\Gamma^{r,y}_{\a,1}}
\frac{dx}{x}\;\I^\a(r,y| x)\;f(x)
\label{eq:def-Ia}
\eeq
with $r$ considered as a parameter.

The following properties of the operator $I^\a$ can be obtained mostly 
from those of $M_{\a\b}^r$ enlisted in the Proposition \ref{th:properties-of-Mab} 
using the correspondence \Ref{eq:M-I}.

\begin{prop}
\hfill\par
\begin{enumerate}
\item The operator $I^\a$ possesses the group property with
respect to the parameter $\a$
\beq
 I^{\a+\b}=I^\a\circ I^\b\,.
\label{eq:I-group}
\eeq

\item
Define an analog $\psi_r^\nu$ of the power function by the formula
(\ref{psi}). Then
\beq
 I^\a:\psi^\nu_r(x)\rightarrow \psi_r^{\nu+\a}(y)\,.
\label{eq:act-Ia-power}
\eeq

\item
For $\a=-g$, $g\in\Z_{\geq0}$ 
the operator $I^\a$ turns into the finite-difference operator
\beq
 I^{-g}:f(x)\rightarrow\sum_{k=0}^g\zeta_{g,k}(r,y)\;f(q^{k-\frac{g}{2}}y)\,,
\label{eq:I^(-g)}
\eeq
\beq
 \zeta_{g,k}(r,y)=(-1)^kq^{-\frac{k(k-1)}{2}}
\left[\matrix{g \cr k}\right]_q
\frac{y^{-2k}\,(1-q^{g-2k}y^{-2})}{(1-q)^g\,(q^{-k}y^{-2};q)_{g+1}}\,
\frac{\Lq{q^{\half}}{r}{y}}{\Lq{q^{\half}}{r}{q^{k-\frac{g}{2}}y}}\,.
\label{eq:def-zeta}
\eeq
Moreover
\beq
 I^{-g}=(I^{-1})^g\,,
\label{eq:(1/I)^g}
\eeq
where
\beq
 I^{-1}:f(x)\rightarrow
\frac{\Lq{q^{\half}}{r}{y}}{(1-q)(1-y^2)}
\left(\frac{f(q^{\half}y)}{\Lq{q^{\half}}{r}{q^{\half}y}}-
\frac{y^2f(q^{-\half}y)}{\Lq{q^{\half}}{r}{q^{-\half}y}}\right).
\label{eq:defI1}
\eeq
\item
In the limit $q\uparrow1$ we have
\beq
 I^{-1}\rightarrow -\frac{y^2}{1-y^2}\frac{d}{dy}\,,
\label{eq:I:q->1}
\eeq
\beq
 \psi_r^\nu(y)\rightarrow \frac{(y+y^{-1}-r-r^{-1})^\nu}{\Gamma(\nu+1)}\,,
\eeq
so, for instance, equality \Ref{eq:act-Ia-power} for $\a=-1$ gives 
the formula
\beq
 I^{-1}\psi_r^\nu=\psi_r^{\nu-1}\,.
\eeq
\end{enumerate}
\label{th:prop-I}
\end{prop}
The identity \Ref{eq:I-group}, like the equivalent identity 
$M_{\a\b}^r:1\rightarrow1$, is a disguised form of the Askey-Wilson (AW)
integral
\Ref{eq:AW-integral}. The formula \Ref{eq:act-Ia-power} is also proved
quite directly using AW integral. The formulas \Ref{eq:I^(-g)} and
\Ref{eq:def-zeta} are equivalent, respectively, to \Ref{eq:M_(-gb)} and
\Ref{def-xi-k}. The equality \Ref{eq:(1/I)^g} is easily proved by 
induction.

The properties of the operator $I^\a$ described in the Proposition 
\ref{th:prop-I} allow to consider it as a fractional power of
the first order difference operator $I^{-1}$ \Ref{eq:defI1}.
In the limit $q\uparrow1$ the operator $I^{-1}$ \Ref{eq:I:q->1} is 
reduced to the pure  differentiation operator by the change of variable
$\xi=y+y^{-1}-2$. The operator $I^\a$, respectively, is reduced to the 
well-known Riemann-Liouville-Weyl fractional integral 
\beq
  I^\a:f(x)\rightarrow\intl_x^\infty dy \;\frac{(y-x)^{\a-1}}{\G(\a)}\;f(y)\,.
\label{eq:Weil-int}
\eeq
In the general case $q\neq1$, however, there is no change of variable
reducing $I^{-1}$ \Ref{eq:defI1} to the pure $q$-derivative
$({\cal D}_qf)(x)=(f(x)-f(qx))/((1-q)x)$, and our operator $I^\a$ \Ref{eq:def-Ia}
cannot be reduced, respectively, to the fractional $q$-integration operator 
studied in \cite{Ag69,AlS66,AV75}.
It is not our task to develop here the complete theory of the operator
$I^\a$ \Ref{eq:def-Ia} although we believe that 
this new $q$-fractional integration operator deserves more comprehensive study. 

\newsection{Separating operator $M_\xi$}
Let us specify the values of the arguments and
parameters in the kernel $\M_{\a\b}(r,y|x)$ \Ref{eq:ker-Mab}
to be
\beq
 \a=\b=g\,, \qquad y=y_-\,, \qquad x=x_-\,, \qquad r=t^{-1}y_+\,, 
\qquad y_\pm=(y_1y_2^{\pm1})^{\half}\,,
\eeq
(in contrast to $\a=g$, $\b=2g$, as for $A_2$ case
\cite{KS2}), and denote the resulting kernel as $\M(y_+,y_-| x_-)$:
\beq
 \M(y_+,y_-| x_-)
=\frac{(1-q)(q;q)_\infty^2\;(x_-^2,x_-^{-2};q)_\infty\;\Lq{t}{y_-}{t^{-1}y_+}}%
{2B_q(g,g)\;\Lq{t^{\half}}{y_-}{x_-}\;\Lq{t^{\half}}{x_-}{t^{-1}y_+}}\,.
\label{4.2}\eeq
Assuming $\xi$ to be an arbitrary complex parameter, we introduce the 
integral operator
$M_\xi$ acting on functions $f(x_1,x_2)$ by the formula
\beq
 (M_\xi f)(y_1,y_2)\equiv\frac{1}{2\pi i}
\int\limits_{\Gamma^{t^{-1}y_+,y_-}_{g,g}}\frac{dz}{z}\;
\M(y_+,y_-|z)\;f(t^{-\half}\xi y_+z,t^{-\half}\xi y_+z^{-1})\,.
\label{4.4}
\eeq
It is apparent from \Ref{4.4} that $M_\xi$ acts on $x_+$ as a simple 
scaling:
\beq
 M_\xi:\phi(x_1x_2)\rightarrow \phi(t^{-1}\xi^2 y_1y_2)\,.
\eeq

Now we can formulate the main result of the paper.

\begin{theo}
The operator $M_\xi$ \Ref{4.4},\Ref{4.2} 
performs the factorisation of (or, in other words, the SoV for) the 
$A_1$ Macdonald polynomials \Ref{eq:def-Pl}:
\beq
 M_\xi:P_\l(x_1,x_2)
\rightarrow F_\l(y_1,y_2)\equiv c_{\l,\xi}\;f_\l(y_1)\;f_\l(y_2)\,,
\label{eq:M-sep-P}
\eeq
where the factorized (or separated) polynomial $f_\l(y)$ is given
in \Ref{eq:def-fl} and the normalization coefficient $c_{\l,\xi}$ is
\beq
 c_{\l,\xi}=t^{-2\l_1+\l_2}\xi^{\abs{\l}}\;\frac{(t;q)_{\l_{21}}}
{(t^2;q)_{\l_{21}}}\,.
\label{eq:def-cl}
\eeq
\label{th:main}
\end{theo}

Note that the product formula \Ref{2.8} for the $q$-ultraspherical polynomials
is equivalent to the formula \Ref{eq:M-sep-P}. To observe it, it is
sufficient to make substitutions $t=\beta$, $\l_{21}=n$, $z=\e^{i\psi}$,
$y_1=t\e^{-2i\theta}$, $y_2=t\e^{-2i\oldphi}$ and to replace the integral over
the unit circle by the integral from $-1$ to $1$.

Taking into account the inversion formulas \Ref{3.9} and \Ref{3.10} for
$M_{\a\b}^r$ we notice that
the kernel of the inverse operator $M_\xi^{-1}$ is obtained
from $\M_{\a\b}(r,y|x)$ \Ref{eq:ker-Mab} by the substitution
$\a=-g$, $\b=2g$, $x:=y_-$, $y:=x_-$, $r:=t^{-\half}\xi^{-1}$
\beq
 \tilde\M(x_+,x_-|y_-)=\frac{(1-q)(q;q)_\infty^2\;(y_-^2,y_-^{-2};q)_\infty
\;\Lq{t^{\half}}{x_-}{t^{-\half}\xi^{-1}x_+}}%
{2B_q(-g,2g)\;\Lq{t^{-\half}}{y_-}{x_-}\;\Lq{t}{y_-}
{t^{-\half}\xi^{-1}x_+}}\,.
\eeq

As the immediate corollary the above remark and the Theorem \ref{th:main} 
we obtain the following result.

\begin{theo}
The inversion of the operator $M_\xi$ \Ref{4.4} is given by the formula
\beq
 (M_\xi^{-1}f)(x_1,x_2)=\frac{1}{2\pi i}
\int\limits_{\Gamma_{-g,2g}^{t^{\tiny -1/2}\xi^{\tiny -1}x_+,x_-}}\frac{dz}{z}\;
\tilde\M(x_+,x_-|z)\;f(t^{\half}\xi^{-1}x_+z,t^{\half}\xi^{-1}x_+z^{-1})\,.
\label{eq:inverse-main}
\eeq
The operator $M_\xi^{-1}$ provides an integral representation for the $A_1$
Macdonald polynomials in terms of the factorized polynomials $f_\l(y)$
\beq
  M^{-1}_\xi:c_{\l,\xi}\;f_\l(y_1)\;f_\l(y_2)
\rightarrow P_\l(x_1,x_2)\,.
\label{eq:int-rel}\eeq
\label{th:inverse-main}
\end{theo}

In contrast to the formula \Ref{eq:M-sep-P} which paraphrases an already known
result, the formula \Ref{eq:int-rel} leads to a new integral relation
for the $q$-ultraspherical polynomials. Note that for positive integer 
$g$ the operator
$M^{-1}_\xi$ becomes a $q$-difference operator (cf. Proposition 
\ref{th:properties-of-Mab}:5).

The proof of the Theorem \ref{th:main} 
mimics that in \cite{KS2} for the $A_2$ case and consists of the
following steps.

\begin{prop}
The operator $M_\xi$ maps bijectively $\Symx$ onto $\Symy$.
\label{th:M-poly-poly}
\end{prop}

\begin{prop}\label{th:q-char-eq-M}
The following operator identity (``quantum characteristic equation'',
cf. \cite{KS2}) is true
 \beq
 (1-qy_j)\,T_{q^2,y_j}\,M_\xi-t^{\half}(1-qt^{-1}y_j)\,T_{q,y_j}\,M_\xi\,H_1
+t(1-qt^{-2}y_j)\,M_\xi\,H_2=0\,,
\label{eq:q-char-eq-M}\eeq
where $H_1$ and $H_2$ are quantum integrals of motion given 
in \Ref{eq:def-H}. 
\end{prop}

\begin{prop}
The separated polynomial  $f_\l(y)$ \Ref{eq:def-fl} is the only
Laurent polynomial solution to the separated equation \Ref{eq:separ-eq}.
\label{th:f-unique}
\end{prop}

\begin{prop}
The normalization coefficient $c_{\l,\xi}$ in \Ref{eq:M-sep-P} is given by
the formula \Ref{eq:def-cl}.
\label{th:calc-cl}
\end{prop}

The main statement of the Theorem \ref{th:main} 
follows then by a standard argument 
(cf. \cite{KS2}). The Propositions \ref{th:M-poly-poly} and 
\ref{th:q-char-eq-M} imply that $M_\xi P_\l$ is a
Laurent polynomial in $y_1$ and $y_2$ satisfying the separated equation
\Ref{eq:separ-eq} in both variables $y_{1,2}$. Then, by virtue
of the Proposition \ref{th:f-unique}, 
it is factorized into $f_\l(y_1)\,f_\l(y_2)$.
It remains only to calculate the normalization coefficient $c_{\l,\xi}$
which is done by the Proposition \ref{th:calc-cl}.

To prove the Proposition \ref{th:M-poly-poly} we shall calculate explicitly
the matrix of $M_\xi$ in special bases in $\Symx$ and $\Symy$.
Define the bases $p_\nu$ and $r_\nu$ in $\Symx$
and, respectively,  $\tilde p_\nu$ and $\tilde r_\nu$ 
in $\Sym(y_1,$ $y_2)$ for
$\nu=\{\nu_1\leq\nu_2\}\in\Z^2$:
\beq
 p_\nu=(x_1x_2)^{\nu_1}(\xi^{-1}x_1,\xi^{-1}x_2;q)_{\nu_{21}}\,, \qquad
 r_\nu=(x_1x_2)^{\nu_2}(t\xi x_1^{-1},t\xi x_2^{-1};q)_{\nu_{21}}\,, 
\label{eq:def-pr}
\eeq
\beq
\tilde p_\nu=(y_1y_2)^{\nu_1}(y_1,y_2;q)_{\nu_{21}}\,, \qquad
 \tilde r_\nu=(y_1y_2)^{\nu_2}(t^2y_1^{-1},t^2y^{-1}_2;q)_{\nu_{21}}\,.
\label{eq:def-tpr}
\eeq
Note that 
\beq
 p_\nu=(x_1x_2)^{\nu_1}R_{00\nu_{21}0}^{g,g}(t^{-1/2}\xi^{-1}x_+,x_-)\,, \qquad
 r_\nu=(x_1x_2)^{\nu_2}R_{000\nu_{21}}^{g,g}(t^{-1/2}\xi^{-1}x_+,x_-)
\eeq
where  $R^{\a\b}_{j_1j_2k_1k_2}$ is defined by \Ref{eq:def-R}. Then,
using \Ref{eq:actMR}, we obtain explicit formulas for the action of $M_\xi$ on
the bases $p_\nu$, $r_\nu$:
\beq
 M_\xi:p_\nu\rightarrow \mu_\nu^{(p)}\;\tilde p_\nu\,,
\qquad
 M_\xi:r_\nu\rightarrow \mu_\nu^{(r)}\;\tilde r_\nu
\label{eq:act-M-p}\eeq
where
\beq
 \mu_\nu^{(p)}=t^{-\nu_1}\xi^{2\nu_1}
\;\frac{(t;q)_{\nu_{21}}}{(t^2;q)_{\nu_{21}}}\,,
\qquad
 \mu_\nu^{(r)}=t^{-\nu_2}\xi^{2\nu_2}
\;\frac{(t;q)_{\nu_{21}}}{(t^2;q)_{\nu_{21}}}\,.
\label{eq:def-mu}
\eeq
The invertibility of $M_\xi$ being obvious from \Ref{eq:act-M-p},
the Proposition \ref{th:M-poly-poly} is thus proved.
Actually, we have obtained much stronger result, having found an effective way
of calculating the action of $M_\xi$ on any polynomial.
\vskip 0.2cm

\noindent
\begin{rem}
There are involutions $U_\xi$ and $V$ in $\Sym(x_1,x_2)$ and $\Sym(y_1,y_2)$, 
respectively,
\beq
  U_\xi:\phi(x_1,x_2)\rightarrow\phi(t\xi^2x_1^{-1},t\xi^2x_2^{-1})\,, \qquad
  V:\phi(y_1,y_2)\rightarrow\phi(t^2y_1^{-1},t^2y_2^{-1})\,,
\label{eq:def-UV}
\eeq
\beq
 U_\xi \;p_\nu=t^{2\nu_1}\xi^{4\nu_1}\;r_{\bar\nu}\,, \qquad
 V\;\tilde p_\nu=t^{4\nu_1}\;\tilde r_{\bar\nu}\,,\qquad
\bar\nu\equiv\{-\nu_2,-\nu_1\}\,,
\label{eq:act-UV-p}
\eeq
which are intertwined by the operator $M_\xi$:
\beq
 M_\xi \;U_\xi=V\;M_\xi\,.
\eeq
\end{rem}

\noindent
\begin{rem} \label{rem:ident-pr}
One can identify $\Symx$ and $\Symy$ setting 
$y_j\equiv\xi^{-1}x_j$. Then $p_\nu\equiv\tilde p_\nu$ becomes the eigenbasis
for $M_\xi$. Similarly, setting $y_j\equiv t\xi^{-1}x_j$ one obtains
$r_\nu\equiv\tilde r_\nu$ as the eigenbasis for $M_\xi$.
\end{rem}
\vskip 0.2cm

\noindent
For  a more detailed discussion of the bases $p_\nu$, $r_\nu$ see the next 
Section.

Now we shall prove the Proposition \ref{th:q-char-eq-M}.
\vskip 0.2cm

\noindent
{\bf Proof.} 
In \cite{KS2} the proof of the identity analogous to 
\Ref{eq:q-char-eq-M} was based on shift relations for the kernel $\M$.
Here we choose another strategy and prove \Ref{eq:q-char-eq-M} purely
algebraically.

Note first that the Hamiltonians $H_j$ \Ref{eq:def-H} are 
described by Jacobian matrices in the  basis $r_\nu$ (not true for $p_\nu$!)
\beq
 H_j\;r_\nu
=a_j^\nu \;r_\nu+b_j^\nu \;r_{\nu_1+1,\nu_2}+c_j^\nu \;r_{\nu_1,\nu_2-1}\,,
\qquad j=1,2\,,
\label{eq:act-H-r}\eeq
where
\beq
\begin{array}{rcl}
 a_1^\nu&=&t^{-\half}\,q^{\nu_1}+t^\half\,q^{\nu_2}\,, \\
 b_1^\nu&=&-t^{-\half}q^{\nu_1}\,(1-q^{\nu_{21}}) \,,\\
 c_1^\nu&=&t^{\frac52}\xi^2q^{-\nu_1+2\nu_2-2}\,(1-q^{\nu_{21}})\,,
\end{array} \qquad
\begin{array}{rcl}
 a_2^\nu&=& q^{\abs{\nu}}\,,\\
 b_2^\nu&=& -q^{\abs{\nu}}\,(1-q^{\nu_{21}})\,,\\
 c_2^\nu&=& t^2\xi^2q^{2\nu_2-2}\,(1-q^{\nu_{21}})\,.
\end{array}
\eeq
Now apply the left-hand-side of \Ref{eq:q-char-eq-M} to the basis $r_\nu$. 
Using the relations 
\begin{subeqnarray}
 \tilde r_{\nu_1+1,\nu_2}&=&\tilde r_\nu\;
(1-q^{\nu_{21}-1}t^2 y_1^{-1})^{-1}(1-q^{\nu_{21}-1}t^2 y_2^{-1})^{-1}, \\
 \tilde r_{\nu_1,\nu_2-1}&=&\tilde r_\nu\; y_1^{-1}y_2^{-1}
(1-q^{\nu_{21}-1}t^2 y_1^{-1})^{-1}(1-q^{\nu_{21}-1}t^2 y_2^{-1})^{-1},
\end{subeqnarray}
\beq
 T_{q,y_j}\;\tilde r_\nu=
\frac{q^{\nu_2}(1-q^{-1}t^2 y_j^{-1})}{1-q^{\nu_{21}-1}t^2 y_j^{-1}}
\;\tilde r_\nu
\eeq
one can express all the $r$'s in terms of $r_\nu$. Discarding the common 
multiplier $r_\nu$ we obtain a rather clumsy rational expression in  $y_1$, 
$y_2$, $t$, $\xi$, $q$ and $q^{\nu_{21}}$ 
which, nevertheless, can be shown to be identically 
zero by a direct computation.
\endproof

The Proposition \ref{th:f-unique} is a direct corollary of the Proposition 12 
from \cite{KS2}.

To complete the proof of the Theorem \ref{th:main} it remains to prove the
Proposition \ref{th:calc-cl}. The proof repeats, with appropriate adjustments,
the proof of the Proposition 7 from \cite{KS2}.

\noindent
{\bf Proof.} Using the variables $x_\pm$ and comparing the asymptotics at
$x_-\rightarrow\infty$ of the monomial symmetric function $m_\l$ 
\Ref{eq:def-monomial}
\beq
 m_\l\sim x_+^{\abs{\l}}\,x_-^{\l_{21}}, \qquad x_-\rightarrow\infty
\eeq
and of the polynomial $r_\nu$ \Ref{eq:def-pr}
\begin{eqnarray}
 r_\nu&\equiv&
x_+^{2\nu_2}\,(t\xi x_+^{-1}x_-^{-1},t\xi x_+^{-1}x_-;q)_{\nu_{21}} \nonumber \\
  &\sim& (-1)^{\nu_{21}}q^{\half\nu_{21}(\nu_{21}-1)}t^{\nu_{21}}\xi^{\nu_{21}}\,
x_+^{\abs{\nu}}\,x_-^{\nu_{21}}, \qquad x_-\rightarrow\infty
\end{eqnarray}
we conclude that the expansion of $m_\l$ in $r_\nu$ is triangular
\beq
 m_\l=(-1)^{\l_{21}}q^{-\half\l_{21}(\l_{21}-1)}t^{-\l_{21}}\xi^{-\l_{21}}
\,r_\l+\sum_{\scriptstyle \nu\prec\l \atop \scriptstyle \nu\neq\l}
(\cdots)\,r_\nu\,,
\label{eq:expand-mr}
\eeq
and, consequently, the expansion of $P_\l$ in $r_\nu$ has the same structure.
Then using \Ref{eq:act-M-p}--\Ref{eq:def-mu} and the asymptotics of 
$\tilde r_\nu$
\beq
 \tilde r_\nu\sim 
(-1)^{\nu_{21}}q^{\half\nu_{21}(\nu_{21}-1)}t^{2\nu_{21}}
\,y_+^{\abs{\nu}}\,y_-^{\nu_{21}}, \qquad y_-\rightarrow\infty
\eeq
we obtain
\beq
 M_\xi\, P_\l\sim t^{-\l_{21}}\xi^{\abs{\l}}\,
\frac{(t;q)_{\l_{21}}}{(t^2;q)_{\l_{21}}}\,y_+^{\abs{\l}}\,y_-^{\l_{21}}\,,
\qquad y_-\rightarrow\infty\,.
\label{eq:asymp-MP}
\eeq
On the other hand, the right hand side of \Ref{eq:M-sep-P} has the asymptotics
\beq
 c_{\l,\xi}\,f_\l(y_1)\,f_\l(y_2)\sim c_{\l,\xi}\,\chi_{\l_1}\,\chi_{\l_2}\,
y_+^{\abs{\l}}\,y_-^{\l_{21}}\,, \qquad y_-\rightarrow\infty\,. 
\label{eq:asymp-cff}
\eeq
Comparing \Ref{eq:asymp-MP} and \Ref{eq:asymp-cff}, and using formula
\Ref{eq:chi-l12} we obtain the desired expression \Ref{eq:def-cl} for
$c_{\l,\xi}$.
\endproof

\newsection{Bases $p_\nu$ and $r_\nu$}
In this section we establish some more properties of the polynomial bases 
$p_\nu$, $r_\nu$, $\tilde p_\nu$, $\tilde r_\nu$ introduced by
\Ref{eq:def-pr}, \Ref{eq:def-tpr} in the previous Section.
Being the eigenfunctions of the separating integral operator $M_\xi$
(see Remark \ref{rem:ident-pr}), these bases must play an important role in 
SoV.

Our first observation is that the bases can be interpreted also
as eigenfunctions of certain quantum integrable systems. The following
statement is easily verified by a direct calculation.

\begin{prop}
 The Laurent polynomials $p_\nu$, $\tilde p_\nu$, $r_\nu$, $\tilde r_\nu$ 
are the eigenfunctions, respectively, of the commuting pairs of
$q$-difference operators
$N_j$, $\tilde N_j$, $Q_j$, $\tilde Q_j$ (j=1,2)
\beq
 N_j\;p_\nu=q^{\nu_j}\;p_\nu\,, 
\quad  \tilde N_j\;\tilde p_\nu=q^{\nu_j}\;\tilde p_\nu\,, \quad
 Q_j\;r_\nu=q^{-\nu_j}\;r_\nu\,, 
\quad  \tilde Q_j\;\tilde r_\nu=q^{-\nu_j}\;\tilde r_\nu\,, 
\eeq
\begin{subeqnarray}
&&\kern-8mm N_1=-\sfrac{x_2(\xi-x_1)}{\xi(x_1-x_2)}\;T_{q,x_1}+
\sfrac{x_1(\xi-x_2)}{\xi(x_1-x_2)}\;T_{q,x_2}\,, \quad
 N_2=-\sfrac{\xi-x_1}{x_1-x_2}\;T_{q,x_1}+\sfrac{\xi-x_2}{x_1-x_2}\;T_{q,x_2}\,,
\\
&&\kern-8mm  \tilde N_1=-\sfrac{y_2(1-y_1)}{y_1-y_2}\;T_{q,y_1}
+\sfrac{y_1(1-y_2)}{y_1-y_2}\;T_{q,y_2}\,, 
\quad
 \tilde N_2=-\sfrac{1-y_1}{y_1-y_2}\;T_{q,y_1}+\sfrac{1-y_2}{y_1-y_2}\;T_{q,y_2}\,,
\\
&&\kern-8mm  Q_1=-\sfrac{x_2(t\xi-x_1)}{t\xi(x_1-x_2)}\;T_{q,x_1}^{-1}+
\sfrac{x_1(t\xi-x_2)}{t\xi(x_1-x_2)}\;T_{q,x_2}^{-1}\,, \quad
 Q_2=-\sfrac{t\xi-x_1}{x_1-x_2}\;T_{q,x_1}^{-1}
+\sfrac{t\xi-x_2}{x_1-x_2}\;T_{q,x_2}^{-1}\,,
\\
&&\kern-8mm  \tilde Q_1=-\sfrac{y_2(t^2-y_1)}{t^2(y_1-y_2)}\;T_{q,y_1}^{-1}
+\sfrac{y_1(t^2-y_2)}{t^2(y_1-y_2)}\;T_{q,y_2}^{-1}\,, 
\quad
 \tilde Q_2=-\sfrac{t^2-y_1}{y_1-y_2}\;T_{q,y_1}^{-1}
+\sfrac{t^2-y_2}{y_1-y_2}\;T_{q,y_2}^{-1}\,.
\end{subeqnarray}
\end{prop}

\noindent
Notice that the above spectral problems are already separated in the variables
$x_j$ ($y_j$).
The integral operator $M_\xi$ obviously intertwines these pairs of first 
order dif\-fe\-rence operators
\beq
 M_\xi \;N_j=\tilde N_j\;M_\xi\,, \qquad  M_\xi \;Q_j=\tilde Q_j\;M_\xi\,.
\label{eq:M-Intrtw-N}
\eeq
Note also that, by virtue of \Ref{eq:def-UV},
\beq
 U_\xi\;N_j\;U_\xi^{-1}=Q_{3-j}\,, \qquad 
 V\;\tilde N_j\;V^{-1}=\tilde Q_{3-j}\,.
\eeq
\vskip 0.2cm

\noindent
\begin{rem}
Choosing the identification $y_j\equiv\xi^{-1}x_j$ of the spaces $\Symx$ and 
$\Symy$, see Remark \ref{rem:ident-pr},
one can identify $N_j\equiv\tilde N_j$. By virtue of \Ref{eq:M-Intrtw-N}
the operator $M_\xi$ commutes then with $N_j\equiv\tilde N_j$ and can be
considered as the time-shift operator for the {\it discrete-time quantum integrable
system} with integrals of motion $N_1$ and $N_2$. 
Same is true for another identification $y_j\equiv t\xi^{-1}x_j$ and the 
operators $Q_j\equiv\tilde Q_j$.
\end{rem}
\vskip 0.2 cm

Having defined several polynomial bases in $\Sym(x_1,x_2)$, such as $P_\l$,
$p_\nu$, $r_\nu$, it is natural to calculate the transition matrices between 
them. The solution is given by the following Proposition. 

\begin{prop}
The expansion of $P_\l$ \Ref{eq:def-Pl} in the bases $p_\nu$ ($r_\nu$) 
\Ref{eq:def-pr} is given by
\beq
  P_\l=\sum_{\nu\prec\l}\pi_\l^\nu\;p_\nu=\sum_{\nu\prec\l}\rho_\l^\nu\;r_\nu
\label{eq:exp-P-pr}
\eeq
where
\begin{subeqnarray}
 \pi_\l^\nu&=&
\frac{(-1)^{\nu_{21}}\xi^{\abs{\l}-2\nu_1}
q^{\frac{\nu_{21}}{2}(\abs{\nu}-2\l_2+1)}
(t;q)_{\l_2-\nu_1}(t;q)_{\nu_2-\l_1}(q;q)_{\l_{21}}}%
{(q;q)_{\l_2-\nu_2}(q;q)_{\nu_1-\l_1}(t;q)_{\nu_{21}}(t;q)_{\l_{21}}
(q;q)_{\nu_{21}}}\,,
\\
 \rho_\l^\nu&=&
\frac{(-1)^{\nu_{21}}(t\xi)^{-2\nu_2+|\l|}q^{\frac{\nu_{21}}{2}(2\l_1+1-|\nu|)}
(t;q)_{\nu_2-\l_1}(t;q)_{\l_2-\nu_1}(q;q)_{\l_{21}}}
{(q;q)_{\l_2-\nu_2}(q;q)_{\nu_1-\l_1}(t;q)_{\nu_{21}}(t;q)_{\l_{21}}
(q;q)_{\nu_{21}}}\,.
\label{eq:formula-pi-rho}
\end{subeqnarray}
The dual expansions of the bases $p_\l$ and $r_\l$ \Ref{eq:def-pr}   
in terms of the $A_1$ Macdonald polynomials $P_\nu$ \Ref{eq:def-Pl} 
are as follows:
\beq
  p_\l=\sum_{\nu\prec\l}Q_\l^\nu\;P_\nu\,,\qquad 
  r_\l=\sum_{\nu\prec\l}R_\l^\nu\;P_\nu
\label{eq:exp-P-pr-qq}
\eeq
where
\begin{subeqnarray}
 Q_\l^\nu&=&
\frac{(-1)^{\nu_{21}}\xi^{2\l_1-\abs{\nu}}q^{\l_1^2-(|\nu|-1)\l_1-\frac{|\nu|}{2}
+\frac{\nu_1^2+\nu_2^2}{2}} (tq;q)_{\l_{21}}(tq;q)_{\nu_{21}}(q;q)_{\l_{21}}}
{(q;q)_{\l_2-\nu_2}(q;q)_{\nu_1-\l_1}(tq;q)_{\nu_2-\l_1}(tq;q)_{\l_2-\nu_{1}}
(q;q)_{\nu_{21}}}\,,\qquad
\\
 R_\l^\nu&=&
\frac{(-1)^{\nu_{21}}(t\xi)^{2\l_2-|\nu|}q^{\l_2^2-(|\nu|+1)\l_2+\frac{|\nu|}{2}
+\frac{\nu_1^2+\nu_2^2}{2}} (tq;q)_{\l_{21}}(tq;q)_{\nu_{21}}(q;q)_{\l_{21}}}
{(q;q)_{\l_2-\nu_2}(q;q)_{\nu_1-\l_1}(tq;q)_{\nu_2-\l_1}(tq;q)_{\l_2-\nu_{1}}
(q;q)_{\nu_{21}}}\,.
\label{eq:formula-RQ}
\end{subeqnarray}
\end{prop}
\Proof
It was already demonstrated in Proposition \ref{th:calc-cl} that the expansion
of $P_\l$ in $r_\nu$ (or $p_\nu$) is triangular with respect to the ordering 
$\prec$ \Ref{eq:order}.
Applying the operators $H_j$ to the equality 
$P_\l-\sum_{\nu\prec\l}\rho_\l^\nu \,r_\nu=0$, using \Ref{eq:act-H-r} and
equating to zero the coefficient at $r_\nu$ we obtain the set of two 
equations
\beq
 \rho_\l^\nu\,(a_j^\nu-h_{\l;j})+\rho_\l^{\nu_1-1,\nu_2}\,b_j^{\nu_1-1,\nu_2}
+\rho_\l^{\nu_1,\nu_2+1}\,c_j^{\nu_1,\nu_2+1}=0\,,\qquad j=1,2\,.
\label{eq:rabc}
\eeq
Solving the above equations with respect to $\rho_\l^{\nu_1-1,\nu_2}$ and
$\rho_\l^{\nu_1,\nu_2+1}$ one obtains the
recurrence relations for the coefficients $\rho_\l^\nu$ 
\begin{subeqnarray}
 \rho_\l^{\nu_1-1,\nu_2}&=&
-\frac{(1-q^{\nu_1-\l_1})(1-tq^{\l_2-\nu_1})}%
{q^{\nu_1-\l_1-1}(1-q^{\nu_{21}+1})(1-tq^{\nu_{21}})}\;\rho_\l^\nu\,,
\\
 \rho_\l^{\nu_1,\nu_2+1}&=&
-\frac{(1-q^{\l_2-\nu_2})(1-tq^{\nu_2-\l_1})}%
{q^{\nu_2-\l_1}t^2\xi^2
(1-q^{\nu_{21}+1})(1-tq^{\nu_{21}})}\;\rho_\l^\nu\,.
\label{eq:recur-rho}
\end{subeqnarray}
As the initial value for the recurrence relations we can use $\rho_\l^\l$
which is extracted directly from the expansion \Ref{eq:expand-mr}:
\beq
\rho_\l^\l=
 (-1)^{\l_{21}}q^{-\half\l_{21}(\l_{21}-1)}t^{-\l_{21}}\xi^{-\l_{21}}\,.
\label{eq:init-cond-rho}\eeq
It is easy to see then that the recurrence relations \Ref{eq:recur-rho}
together with the initial condition \Ref{eq:init-cond-rho} possess the unique
solution (\ref{eq:formula-pi-rho}b) for $\nu\prec\l$. Note that the solution
is compatible with the triangularity condition $\nu\prec\l$ since from
\Ref{eq:recur-rho} it follows that $\rho_\l^\nu=0$ for $\nu_1<\l_1$ or
$\nu_2>\l_2$.

Similarly, for the $R$-coefficients we obtain from 
$r_\l=\sum_{\nu\prec\l}R_\l^\nu \,P_\nu$ the equations
\beq
 (a_j^\l-h_{\l;j})\,R_\l^\nu+b_j^\l \,R_{\l_1+1,\l_2}^\nu
+c_j^\l \,R_{\l_1,\l_2-1}^\nu=0
\label{eq:Rabc}
\eeq
which give rise to the recurrence relations
\begin{subeqnarray}
 R_{\l_1+1,\l_2}^{\nu}&=&
\frac{(1-q^{\nu_1-\l_1})(1-tq^{\nu_2-\l_1})}%
{(1-q^{\l_{21}})(1-tq^{\l_{21}})}\;R_\l^\nu\,,
\\
 R_{\l_1,\l_2-1}^{\nu}&=&
\frac{(1-q^{\l_2-\nu_2})(1-tq^{\l_2-\nu_1})}%
{q^{2\l_2-|\nu|-2}t^2\xi^2
(1-q^{\l_{21}})(1-tq^{\l_{21}})}\;R_\l^\nu
\end{subeqnarray}
which, in turn, being combined with the initial condition
\beq
 R_\l^\l=(-1)^{\l_{21}}q^{\frac12 \l_{21}(\l_{21}-1)}t^{\l_{21}}\xi^{\l_{21}}
\label{eq:init-cond-R}
\eeq
produce the formula (\ref{eq:formula-RQ}b).

The easiest way to obtain the corresponding formulas for $\pi_\l^\nu$ and 
$Q_\l^\nu$ is to use the involution $U_\xi$ \Ref{eq:def-UV}.
Noting that
\beq
 U_\xi\; P_\l=t^{\abs{\l}}\xi^{2\abs{\l}}\;P_{\bar\l}
\eeq
and using the equality $ U_\xi\,p_\nu=t^{2\nu_1}\xi^{4\nu_1}\,r_{\bar\nu}$
\Ref{eq:act-UV-p} one obtains the formulas
\beq
 \pi_\l^\nu=(t\xi^2)^{\abs{\l}-2\nu_1}\rho_{\bar\l}^{\bar\nu}
\eeq
and
\beq
 Q_\l^\nu=(t\xi^2)^{2\l_1-\abs{\nu}}R_{\bar\l}^{\bar\nu}
\eeq
which lead directly to the expressions (\ref{eq:formula-pi-rho}a) and
(\ref{eq:formula-RQ}a). 
\endproof
\vskip 0.2cm

\noindent
{\bf Corollary.} The matrices $\rho_\l^\nu$ and $R_\l^\nu$, resp.\ 
$\pi_\l^\nu$ and $Q_\l^\nu$, being mutually inverse, we obtain 
the following algebraic identities:
\be
  \sum_{\mu\prec\nu\prec\l}R_\l^\nu\;\rho_\nu^\mu=0\,,\quad
  \sum_{\mu\prec\nu\prec\l}\rho_\l^\nu\;R_\nu^\mu=0\,,\quad
  \sum_{\mu\prec\nu\prec\l}Q_\l^\nu\;\pi_\nu^\mu=0\,,\quad
  \sum_{\mu\prec\nu\prec\l}\pi_\l^\nu\;Q_\nu^\mu=0\,,
\label{ha-ha}\ee
if $\mu\prec\l$ and $\mu\neq\l$, otherwise
$R_\l^\l\;\rho_\l^\l=1\,,\;Q_\l^\l\;\pi_\l^\l=1$. 
\vskip 0.2cm

\noindent
The identities can also  be  verified directly
with the help of the $q$-Vandermond formula (\cite{GR}, formula (1.5.3)):
$$
{}_2\oldphi{}_1\left[\begin{array}{c}q^{-n},b \\ c\end{array};q,q\right]=
\frac{(c/b;q)_n}{(c;q)_n}\;b^n\,.
$$

Using the relation \Ref{eq:M-sep-P} together with \Ref{eq:act-M-p}
it is easy to transform the expansions \Ref{eq:exp-P-pr} and \Ref{eq:exp-P-pr-qq}
into, respectively,
\beq
 F_\l=\sum_{\nu\prec\l}\tilde\pi_\l^\nu\,\tilde p_\nu
     =\sum_{\nu\prec\l}\tilde\rho_\l^\nu\,\tilde r_\nu
\label{eq:expand-F}\eeq
and
\beq
 \tilde p_\l=\sum_{\nu\prec\l}\tilde Q_\l^\nu \,F_\nu, \qquad
 \tilde r_\l=\sum_{\nu\prec\l}\tilde R_\l^\nu \,F_\nu
\eeq
where
\beq
 \tilde\pi_\l^\nu=\pi_\l^\nu \,\mu_\nu^{(p)}, \qquad
 \tilde\rho_\l^\nu=\rho_\l^\nu \,\mu_\nu^{(r)}
\eeq
and
\beq
 \tilde Q_\l^\nu=Q_\l^\nu \,(\mu_\l^{(p)})^{-1}, \qquad
 \tilde R_\l^\nu=R_\l^\nu \,(\mu_\l^{(r)})^{-1}
\eeq
(for the definition of $\mu_\l^{(p)}$ and $\mu_\l^{(r)}$ 
see \Ref{eq:def-mu}).

The second expansion in \Ref{eq:expand-F} can also be found 
in the book \cite{GR}. Indeed, consider the formula (i) in the 
Excersize 8.1 from \cite{GR}, i.e. the expansion of the 
product of two little $q$-Jacobi polynomials in the form
$$
p_n(x_1;a,b;q)\;p_n(x_2;a,b;q)=(-aq)^nq^{\left({n \atop 2}\right)}
\;\frac{(bq;q)_n}{(aq;q)_n}\;
\sum_{m=0}^n\frac{(q^{-n},abq^{n+1},x_1^{-1},x_2^{-1};q)_m}
{(q,bq;q)_m}\;\times
$$
$$
\times\;\left(\frac{x_1x_2q^{1-m}}{a}\right)^m
\;\sum_{k=0}^m\frac{(q^{-m},b^{-1}q^{-m};q)_k}
{(q,aq,x_1q^{1-m},x_2q^{1-m};q)_k}\;
(abx_1x_2)^k\;q^{k^2+2k}\,.
$$
If we substitute now $n=\l_{21}$, $x=t^{-2}y$, $a=t^{-1}q^{-\l_{21}}$,
$b=t^2q^{-1}$ and use the formula \Ref{eq:def-fl}
for the separated polynomials $f_\l(y)$ then we arrive to the expansion 
\Ref{eq:expand-F} of the polynomial $F_\l$ in terms of the 
polynomials $\tilde r_\nu$. 

\newsection{Concluding remarks}
We have demonstrated that even in the simplest case of the $A_1$ type Macdonald
polynomials the method of separation of variables   
is capable to provide a new interpretation of the product formula \Ref{2.8}
for classical
orthogonal polynomials of one variable (Theorem \ref{th:main}) and even
to produce a new result: the integral relation \Ref{eq:inverse-main} for the
same polynomials (Theorem \ref{th:inverse-main}).

The SoV approach proposed in \cite{S1,KS1,KS2} and in
the present paper seems to be quite general and
could provide thus a useful tool for studying the 
integral relations between various one- and multidimensional special functions.
The most interesting next application of the method would be the case of the
$A_n$ Macdonald polynomials, $n>2$. The work in this direction is in progress.
We hope also to apply the same method to the pair $(A_2,A_1)$
in the cases of periodic Toda lattice and elliptic Calogero-Moser
and Ruijsenaars models elsewhere, thereby 
getting some (new) integral relations for the 
Mathieu and Lam\'{e} types of functions and for their multivariable 
(and $q$-) analogs. 

\section*{Acknowledgments}
We thank I.V.~Komarov for discussion on the subject
of product formulas and B.D.~Sleeman for discussion 
on various aspects of the multiparameter spectral problem. 
VBK wishes to acknowledge the support of EPSRC. 

\setcounter{section}{0}\renewcommand{\thesection}{\Alph{section}}
\newappendix{Classical Ruijsenaars system}
An integrable Hamiltonian system with $n$ degrees of freedom is determined
by a $2n$-dimensional symplectic manifold (phase space) and $n$ independent 
functions (Hamiltonians) $H_j$ commuting with respect to the Poisson bracket
\beq
 \{H_j,H_k\}=0\,,   \qquad j,k=1,\ldots,n\,.
\eeq
For our purposes it is more convenient to use Weil canonical variables
$(T_{x_j},x_j)$ 
\beq
 \{T_{x_j},T_{x_k}\}=\{x_j,x_k\}=0\,, \quad \{T_{x_j},x_k\}=
-i\,T_{x_j}\,x_k\,\d_{jk}\,, \qquad j,k=1,\ldots,n
\label{eq:def-Weyl-can}
\eeq
rather then the traditional coordinates--momenta.

To find a SoV means then to find
a canonical transformation $M:(x,T_x)\mapsto (y,T_y)$,
$\;M:H_j^{(x)}\mapsto H_j^{(y)}$ such that there exist $n$ relations 
\beq
\Phi(y_j,T_{y_j};H_1^{(y)},\ldots,H_n^{(y)})=0\,,
\qquad j=1,\ldots,n\,,
\label{eq:sep-var-cl}
\eeq
separating the variables $y_j$.
The most common way to describe a canonical transformation is the one in terms 
of its generating function $F(y|x)$.

Presently, no algorithm is known of constructing a SoV for
any given integrable system. Nevertheless, 
there exists a fairly effective practical recipe
based on the classical inverse scattering method.
A detailed description of the procedure with many examples can be found in 
the review paper \cite{S1}. Here we describe very briefly its main steps
and apply it afterwards to the classical trigonometric 
$2$-particle Ruijsenaars system. 

A Lax matrix for a given integrable system
is a matrix $L(u)$ dependent
on a ``spectral parameter'' $u\in \C$ such that its characteristic 
polynomial obeys two conditions
\bea
&&(i) \quad\mbox{{\it Poisson involutivity:}}\nonumber\\
&&\quad\quad \{\det(L(u)-v\cdot\id),\det(L(\tilde u)-\tilde v\cdot\id)\}=0\,,
\quad\forall u,\tilde u,v,\tilde v\in\C; \nonumber\\
&&(ii)\quad \det(L(u)-v\cdot\id) \quad \mbox{{\it generates all integrals
of motion}}\;\;H_i\,.
\nonumber\eea
A Baker-Akhiezer (BA) function is the eigenvector 
\beq
L(u)\;f(u)=v(u)\;f(u)
\label{5.5}\eeq
of the Lax matrix $L(u)$, provided that
a normalization of the eigenvectors $f(u)$ is fixed
\beq
\sum_{i=1}^n\a_i(u)\;f_i(u)=1\,,\qquad (\;f(u)\equiv
(f_1(u),\ldots,f_n(u))^t\;)\,.
\label{5.6}\eeq
Supposing $L(u)$ to be a rational function in $u$,
the pair $(u,v)$ can be thought of as a point of the algebraic 
{\it spectral curve} 
\beq
 \det(L(u)-v\cdot\id)=0\,.
\label{eq:spectral-curve}
\eeq
The BA function $f(u)$ is then a meromorphic function on the spectral curve.

The recipe for finding an SoV
is simple:
\vskip 0.2cm
\noindent
{\it The separation variables $y_j$ are poles of the Baker-Akhiezer 
function, provided it is properly normalized.
The corresponding eigenvalues $T_{y_j}$ of $L(y_j)$, or some functions  
of them, serve as the canonically conjugated variables.
}
\vskip 0.2cm

It is easy to see that the pairs $(y_j,T_{y_j})$ thus defined satisfy the 
separation equations \Ref{eq:sep-var-cl} for 
$\Phi_j\equiv\det(L(y_j)-T_{y_j}\cdot\id)$. The canonicity of the variables 
$(y_j,T_{y_j})$ should be verified independently.
No general recipe is known how to guess the proper (that is producing
canonical variables) normalization for the BA function. 
In many cases the simplest {\it standard} normalization 
$\vec\a=\vec\a_0\equiv(0,\ldots,0,1)$ works.
In other cases the vector $\vec\a$ may depend on the
spectral parameter $u$ and the dynamical variables $(T_{x_j},x_j)$. We shall 
refer to such normalization as a {\it dynamical} one.

To be more concrete, let
$f_i^{(j)}=\res_{u=y_j}f_i(u)$ and $T_{y_j}\equiv v(y_j)$.
Then from (\ref{5.5})--(\ref{5.6}) we have the overdetermined system
of $n+1$ linear homogeneous equations for $n$ components $f_i^{(j)}$
of the vector $f^{(j)}$:
\beq
\left\{\begin{array}{l}
L(y_j)\;f^{(j)}=T_{y_j}\;f^{(j)}\,, \\
\sum_{i=1}^n\a_i(y_j)\;f_i^{(j)}=0\,.\end{array}\right.
\label{5.7}\eeq
The pair $(u,v)\equiv(y_j,T_{y_j})$ is thus determined from the condition
\beq
\mbox{{\rm rank}}\pmatrix{\vec\a(u)\cr L(u)-v\cdot\id}=n-1
\label{5.9}\eeq
where $\vec\a$ is understood as a row-vector.
Finally, the condition (\ref{5.9}) can be rewritten as the 
following vector equation:
\beq
\vec \a\cdot (L(u)-v\cdot\id)^\wedge=0\,,
\label{5.11}\eeq
where the wedge  denotes the matrix of cofactors. 
Eliminating $v$ from (\ref{5.11}) 
we get an algebraic equation for $y_j$
\beq
B(u)=\det\pmatrix{\vec\a\cr
\vec\a \cdot L(u)\cr
\vdots\cr
\vec\a \cdot L^{n-1}(u)}=0\,.
\label{5.12}\eeq
Also, from equations (\ref{5.11}) we can get formulas for $v$ in the form
$$
v=A(u)
$$
with $A(u)$ being some rational functions of the entries of $L(u)$. 

To validate the choice of normalization $\vec\a(u)$
it remains to verify (somehow) the canonicity of brackets between
the whole set of separation variables, namely: between zeros $y_j$ of $B(u)$ 
and their conjugated variables $T_{y_j}\equiv v(y_j)=A(y_j)$. 
To do this final
calculation one needs information about Poisson brackets between
entries of the Lax matrix $L(u)$. 

Now we apply the above scheme to the $2$-particle Ruijsenaars system. 
But, first, let us recall the definition of the $n$-particle system. 

The $n$-particle ($A_{n-1}$ type) trigonometric Ruijsenaars model \cite{Ru}
is formulated in terms of the Weyl canonical system of coordinates
$(T_{x_j},x_j)$, 
$\abs{x_j}=1$, $T_{x_j}\in\R$ $(j=1,2,\ldots,n)$
with the Poisson brackets \Ref{eq:def-Weyl-can}.
The mutually Poisson commuting integrals of motion
$H_i$ are defined as follows
\beq
 H_i=\sum_{J\subset\{1,\ldots,n\} \atop \left|J\right|=i}
  \left(\prod_{j\in J \atop k\in \{1,\ldots,n\}\setminus J} v_{jk} \right)
  \left(\prod_{j\in J}T_{x_j} \right)\,, \qquad i=1,\ldots,n\,,
\label{eq:def-HH}
\eeq
where
\beq
 v_{jk}=\frac{t^{\half}x_j-t^{-\half}x_k}{x_j-x_k}\,, \qquad
 0\leq t\leq 1\,,
\label{eq:def-ww}
\eeq
$t$ being a parameter of the model. 
The $n\times n$ Lax matrix (or the $L$-operator) is given by the formula
\beq
   L(u)=D(u)E(u)
\label{eq:L=DE}
\eeq
where
\beq
 D_{jk}=\frac{(1-t)(t^n-u)}{2t^{\frac{n+1}{2}}(1-u)}
\left(\prod_{i\neq j}v_{ji}\right)T_{x_j}\d_{jk}\,,\qquad
  E_{jk}=\frac{t^n+u}{t^n-u}-\frac{tx_j+x_k}{tx_j-x_k}\,.
\label{eq:def-E-gen}
\eeq
The characteristic polynomial of the matrix $L(u)$ (\ref{eq:L=DE})
generates the integrals (\ref{eq:def-HH})
\begin{eqnarray}
\lefteqn{
 (-1)^nt^{\frac{n(n-1)}{2}}(t^n-u)(1-u)^n\det(z\cdot\id-L(u))
} \nonumber \\
&&=\sum_{k=0}^n(-1)^kt^{\frac{n-1}{2}k}(t^k-u)(1-u)^k(t^n-u)^{n-k}
H_{n-k}z^k
\label{eq:char-poly-N}
\end{eqnarray}
where it is assumed $H_0\equiv1$.

In \cite{KS2} we considered the case $n=3$ and constructed a
classical (as well as quantum) SoV
for it, giving explicitely the generating function of the 
separating canonical transform in terms of Euler's dilogarithm function
\cite{Ki}
\beq
 \dilog(z):=-\int_0^z\frac{dt}{t}\ln(1-t)=\sum_{k=1}^\infty\frac{z^k}{k^2}\,.
\label{eq:dilog}
\eeq
In the present paper we are concerned with the $2$-particle (or the $A_1$)
case. Putting $n=2$ we have, respectively:
\beq
 D(u)=\frac{(1-t)(t^2-u)}{2t^{\frac32}(1-u)}
\pmatrix{v_{12}T_{x_1} & 0 \cr 0 & v_{21}T_{x_2}},\quad
  E_{jk}=\frac{t^2+u}{t^2-u}-\frac{tx_j+x_k}{tx_j-x_k}\,.
\label{eq:def-D(u)}
\eeq
The characteristic polynomial takes the form
\beq
t(1-u)\,\det(z\cdot\id-L(u))
=t(1-u)\,z^2-t^{\frac12}(t-u)\,H_1\,z+(t^2-u)\,H_2
\label{eq:char-poly-cl}
\eeq
where the integrals of motion $H_j$ are given by the formulas
\beq
 H_1=\frac{t^{\half}x_1-t^{-\half}x_2}{x_1-x_2}\;T_{x_1}+
\frac{t^{\half}x_2-t^{-\half}x_1}{x_2-x_1}\;T_{x_2}\,, 
\qquad H_2=T_{x_1}\;T_{x_2}\,.
\label{eq:def-H2}
\eeq

The standard choice of the normalization $\vec\a=\vec\a_0\equiv(0,1)$ leads
to a trivial SoV in the coordinates $x_\pm=(x_1x_2^{\pm1})^{1/2}$ (separating
the center-of-mass motion). There exists, however, another, more complicated
SoV.
Let us choose the following normalization $\vec\a$ of the Baker-Akhiezer 
function 
\beq
 \vec\a\equiv (\a_1,\a_2)=\left(
\frac{t^2+u}{t^2-u}-\frac{t\xi+x_1}{t\xi-x_1}\,,\quad
\frac{t^2+u}{t^2-u}-\frac{t\xi+x_2}{t\xi-x_2}\right),
\label{eq:def-alpha}
\eeq
where $\xi\in\C$ is a free parameter.
Notice that the chosen normalization is dynamical and dependent on 
the spectral parameter $u$. 

Now introduce the functions $A_j(u)$
\beq
 A_1(u)=L_{11}(u)-\frac{\a_1}{\a_2}\;L_{12}(u)\,, \qquad
 A_2(u)=L_{22}(u)-\frac{\a_2}{\a_1}\;L_{21}(u)\,,
\label{eq:def-A}
\eeq
\beq
 A_j(u)=T_{x_j}\;a_j(u)\,,
\eeq
\beq
 a_j(u)=\frac{(t^2-u)(\xi-x_j)(\xi u-x_{3-j})}%
{(1-u)(t\xi-x_j)(\xi u-tx_{3-j})}\,, \qquad j=1,2\,.
\label{eq:def-a}
\eeq
The separation variables $(T_{y_j},y_j)$, $j=1,2$, are defined
then by the equations
\beq
 T_{y_j}=A_1(y_j)=A_2(y_j)\,, \qquad j=1,2\,;
\label{eq:def-y}
\eeq
which can be put, alternatively, as four equations
for four variables $y_1, y_2, T_{y_1}, T_{y_2}$ in the form
\beq
 T_{y_j}=T_{x_k}\;a_k(y_j)\,, \qquad j,k=1,2\,.
\eeq
{}From the invariance of $a_1(u)/a_2(u)$ with respect to 
the change $u\rightarrow tx_1x_2/u\xi^2$ we can deduce
the relation
\beq
 y_1y_2\xi^2=tx_1x_2\,.
\label{eq:constr-yt}
\eeq
Having checked by straightforward calculations the relation
\beq
 t(1-y)\,a_1(y)\,a_2(y)-t^{\frac12}(t-y)\,[\,v_{12}a_2(y)+v_{21}a_1(y)\,]
+(t^2-y)=0
\label{eq:char-eq-a-cl}
\eeq
for $a_j$'s, we obtain at once the separation equation 
for each pair of $(T_y,y)$ variables
\beq
 t(1-y)\;T_y^2-t^{\frac12}(t-y)\;H_1\;T_y+(t^2-y)\;H_2=0\,.
\eeq
The (Weyl type) canonical Poisson brackets \Ref{eq:def-Weyl-can} for the variables
$T_{y_j},y_j$ are easily checked on a computer. 
The other way to establish that the separating transformation from
$(T_{x_j},x_j)$ to $(T_{y_j},y_j)$ is canonical is to give an explicit
generating function for it which 
has the simplest form in terms of new (canonical) $\pm$-variables:
\beq
 x_\pm=x_1^{\frac12}x_2^{\pm\frac12}\,, \qquad
 T_x^\pm=T_{x_1}T_{x_2}^{\pm1}\,,
\eeq
\beq
 y_\pm=y_1^{\frac12}y_2^{\pm\frac12}\,, \qquad
 T_y^\pm=T_{y_1}T_{y_2}^{\pm1}\,,
\eeq
\beq
 x_+=t^{-\frac12}\xi \,y_+\,.
\label{ext}\eeq
The generating function $F(T_{y}^+,y_-|x_+,x_-)$ has the form
\beq
 F(T_{y}^+,y_-|x_+,x_-)=
i\ln T_{y}^+\;\ln(\xi^{-1}t^{\frac12}x_+)+\tilde F(y_+,y_-|x_+,x_-)\,,
\eeq
\bea
 \tilde F&=&
i\left({\cal L}(t^{\half};y_-,x_-)
 +{\cal L}(t^{\half};t^{-\half}\xi^{-1}x_+,x_-)
 -{\cal L}(t;t^{-\half}\xi^{-1}x_+,y_-)\right)\nonumber \\
&& -i\dilog(x_-^2)-i\dilog(x_-^{-2})\,,
\label{aaa}\eea
where we have introduced the notation
\beq
 {\cal L}(\nu;x,y):=\dilog(\nu xy)+\dilog(\nu xy^{-1})
+\dilog(\nu x^{-1}y)+\dilog(\nu x^{-1}y^{-1})\,.
\label{eq:def-LL}
\eeq
Function $\tilde F(y_+,y_-|x_+,x_-)$ (\ref{aaa}) satisfies the equations
\beq
 x_+\dd_{x_+}\tilde F=i\ln\frac{T_{x}^+}{T_{y}^+}\,, \quad
 x_-\dd_{x_-}\tilde F=i\ln T_{x}^-\,, \quad
 y_-\dd_{y_-}\tilde F=-i\ln Y_y^-\,, \quad
 y_+\dd_{y_+}\tilde F=0\,.
\eeq

There is striking similarity between
the classical SoV performed above for the case $n=2$ and 
the one for $n=3$ constructed in \cite{KS2}. Indeed, we recall
that there we had the generating function $F(T_{y}^+,y_-|x_+,x_-)$ 
((2.31) of \cite{KS2}) of the form 
$F:=i\ln T_y^+\ln(t^{\frac32}x_+)+\tilde F$ where
\begin{eqnarray}
 \tilde F&:=&{i}\bigl({\cal L}(t^{\half};y_-,x_-)
+{\cal L}(t;x_+,x_-)-{\cal L}(t^{\frac32};x_+,y_-) \nonumber \\
&&-\dilog(x_-^2)-\dilog(x^{-2}_-)\bigr)\,,
\label{eq:def-F}
\end{eqnarray}
and an extra condition which replaced (\ref{ext}) was as follows:
$$
x_+=t^{-\frac32}y_+\,,
$$
while the definition of the $x_\pm$-variables was a bit different:
$$
x_+=x_1^{1/2}x_2^{1/2}x_3^{-1}, \qquad x_-=x_1^{1/2}x_2^{-1/2}\,.
$$
Actually, there is a simple and an elegant explanation of this 
similarity; moreover, in the general case of $n$ degrees of freedom
we could say that a SoV for the $A_{n-1}$ problem with the standard
normalization vector $\vec\a_0\equiv (0,0,\ldots,0,1)$ 
implies a SoV for the $A_{n-2}$ problem with the non-standard
normalization vector $\vec\a=(\a_1,\ldots,\a_{n-1})$ (cf. (\ref{eq:def-alpha})
for $n=3$),
\beq
\vec\a=\left(
\frac{t^{n-1}+u}{t^{n-1}-u}-\frac{t\xi+x_1}{t\xi-x_1}\,,\ldots\,,
\frac{t^{n-1}+u}{t^{n-1}-u}-\frac{t\xi+x_{n-1}}{t\xi-x_{n-1}}\right),
\label{kk}\eeq
if we choose $\xi=x_n$. Let us demonstrate this explicitly. 

If we remove the last ($n$th) row and the last column from the Lax matrix
$L^{(n)}(u)\equiv L(u)$ (\ref{eq:L=DE})--(\ref{eq:def-E-gen}) then,
as one can check easily,
we get $(n-1)\times (n-1)$ Lax matrix $\tilde L(u)$ for an
integrable system of $n-1$ degrees of freedom which is equivalent
through a simple canonical transformation to the 
standard $(n-1)$-particle Ruijsenaars model.
This 1-degree-of-freedom-less
system with the Lax matrix $L^{(n-1)}(u)$ 
directly inherits the non-standard SoV with the normalization (\ref{kk})
from the standard one (with the $\vec\a_0$) for the system with $n$ degrees
of freedom. Indeed, to see this, it is sufficient to note that 
separation variables $(T_y,y)$
for {\it both} systems (one with $L^{(n)}(u)$ and one with
$\tilde L(u)$) are defined from the intersection
of two spectral curves:
\beq
\left\{\matrix{\det(L^{(n)}(y)-T_y\cdot\id)=0\,,\cr 
\det(\tilde L(y)-T_y\cdot\id)=0\,.}\right.
\label{5kk}\eeq
In other words, the condition of the standard SoV for the first problem,
\beq
\mbox{{\rm rank}}\pmatrix{\vec\a_0\cr L^{(n)}(y)-T_y\cdot\id}=n-1\,,
\label{6kk}\eeq
implies the following condition of SoV for the second problem:
\beq
\mbox{{\rm rank}}\pmatrix{\vec\a_1\cr \tilde L(y)-T_y\cdot\id}=n-2\,,
\label{7kk}\eeq
where 
$$
\vec\a_1=\left(
\frac{t^{n}+u}{t^{n}-u}-\frac{tx_n+x_1}{tx_n-x_1}\,,\ldots\,,
\frac{t^{n}+u}{t^{n}-u}-\frac{tx_n+x_{n-1}}{tx_n-x_{n-1}}\right).
$$
Applying this to our $(A_2,A_1)$ pair,
we have proved in \cite{KS2} the standard
SoV for the $A_2$ problem ($L^{(3)}(u)=L(u),\;n=3$), hence, the same
separation variables can work also for the related $\tilde L(u)$ problem
which, in turn, is related to the initial 2-degrees-of-freedom ($A_1$)
problem (with the Lax matrix $L^{(2)}(u)=L(u),\;n=2$) through the following 
canonical transformation:
\beq
y=t^{-1}\tilde y\,,\qquad T_{x_j}=t^\half v_{j3}\;\tilde T_{x_j}\,,\qquad
\xi=x_3\,,\qquad T_y=t^\half \frac{1-\tilde y}{t^\half-t^{-\half}\tilde y}\;
\tilde T_y\,.
\label{final}\eeq

The analogous transformation connects the corresponding quantum
$A_2$ and $A_1$ systems. Indeed, let us establish an equivalence of the 
separating kernel $\M^{(2)}$ constructed in the Section 4
for the case of the $A_1$ Macdonald polynomials
and the $\M^{(3)}$ from \cite{KS2} for the 3-variable case.
Let us write down both kernels:
for the $2$-variable case 
\beq
 \M ^{(2)}\equiv \M_{g,g}(r,y|x)=
\frac{(1-q)\;(q,q,x^2,x^{-2};q)_\infty\;\Lq tyr}%
{2B_q(g,g)\;\Lq{t^\half}xy\;\Lq{t^\half}xr}
\eeq
and for the $3$-variable case
\beq
 \M ^{(3)}\equiv \M_{g,2g}(\tilde r,y|x)=
\frac{(1-q)\;(q,q,x^2,x^{-2};q)_\infty\;\Lq{t^{\frac32}}y{\tilde r}}%
{2B_q(g,2g)\;\Lq{t^\half}xy\;\Lq tx{\tilde r}}\,.
\eeq
Substituting 
\beq
 y=y_-\,, \qquad x=x_-\,, \qquad
r=t^{-1}y_+=t^{-\half}\xi^{-1}x_+\,,
\eeq
and, respectively,
\beq
 y=\tilde y_-\,, \qquad x=x_-\,, \qquad
\tilde r=t^{-\frac32}\tilde y_+=x_+x_3^{-1}\,,
\eeq
and also
\beq
 B_q(a,b)=(1-q)\frac{(q,q^{a+b};q)_\infty}{(q^a,q^b;q)_\infty}\,,
\eeq
we get
\beq
 \M^{(2)}=\frac{\Lq t{y_-}{t^{-1}y_+}}{(t^2;q)_\infty}\;
\frac{(q,t,x_-^2,x_-^{-2};q)_\infty}{2\Lq{t^\half}{x_-}{y_-}}\;
\frac{(t;q)_\infty}{\Lq{t^\half}{x_-}{t^{-\half}\xi^{-1}x_+}}\,,
\eeq
\beq
 \M^{(3)}=\frac{\Lq{t^{\frac32}}{y_-}{t^{-\frac32}\tilde y_+}}%
{(t^3;q)_\infty}\;
\frac{(q,t,x_-^2,x_-^{-2};q)_\infty}{2\Lq{t^\half}{x_-}{y_-}}\;
\frac{(t^2;q)_\infty}{\Lq t{x_-}{x_+x_3^{-1}}}\,.
\eeq
Identifying the variables
\beq
 \xi\equiv x_3\,, \qquad \tilde y_j\equiv t\,y_j\,, \quad j=1,2\,,
\eeq
we obtain the equivalence of two operators: 
\beq
 M^{(2)}=W_y^{-1}\circ M^{(3)}\circ W_x
\eeq
where
\beq
W_x=\frac{(t;q)_\infty\;\Lq t{x_-}{x_3^{-1}x_+}}%
{(t^2;q)_\infty\;\Lq{t^\half}{x_-}{t^{-\half}\xi^{-1}x_+}}=
\frac{(t,tx_1x_3^{-1},tx_2x_3^{-1};q)_\infty}%
{(t^2,x_1x_3^{-1},x_2x_3^{-1};q)_\infty}\,,
\eeq
\beq
W_y=\frac{(t^2;q)_\infty\;\Lq{t^{\frac32}}{y_-}{t^{-\frac32}\tilde y_+}}%
{(t^3;q)_\infty\;\Lq t{y_-}{t^{-2}\tilde y_+}}=
\frac{(t^2,\tilde y_1,\tilde y_2;q)_\infty}%
{(t^3,t^{-1}\tilde y_1,t^{-1}\tilde y_2;q)_\infty}\,.
\eeq
In particular, the relations between $q$-shift operators are as follows
(cf. (\ref{final})):
\beq
W_x\circ T_{q,x_j}\circ W_x^{-1}=\frac{1-tx_jx_3^{-1}}{1-x_jx_3^{-1}}\;
\tilde T_{q,x_j}=t^\half v_{j3}\;\tilde T_{q,x_j}\,,
\eeq
\beq
W_y\circ T_{q,y_j}\circ W_y^{-1}=\frac{1-\tilde y_j}{1-t^{-1}\tilde y_j}\;
\tilde T_{q,y_j}=t^\half\frac{1-\tilde y_j}{t^\half-t^{-\half}\tilde y_j}
\;\tilde T_{q,y_j}\,.
\eeq
\newpage

\end{document}